\documentclass[12pt,a4paper]{article}
\usepackage{amsmath}
\usepackage{amsfonts}
\usepackage[english,french]{babel}
\usepackage{fancyhdr}
\newcommand{\datum}{Version 99.07.23}
\date{June 1999}
\title{Equity Allocation and Portfolio Selection in Insurance}
\author{Erik Taflin\footnote{AXA,
23, Avenue Matignon, 75008 Paris, France; erik.taflin@u-bourgogne.fr,
erik.taflin@axa.com}}
\newtheorem{theorem}{Theorem}[section]
\newtheorem{lemma}[theorem]{Lemma}

\newtheorem{remark}[theorem]{Remark}
\numberwithin{equation}{section}
\begin{document}
\selectlanguage{english}

%%%%%%%%%%%%%%%%%%%%%%%%%%%%
\maketitle
\thispagestyle{empty}
%%%%%%%%%%%%%%%%%%%%%%%%%%%
\begin{abstract}
\noindent
A discrete time probabilistic model, for optimal equity allocation and portfolio
selection, is formulated so as to apply to (at least) reinsurance. In the
context of a company with several portfolios (or subsidiaries), representing
both liabilities and assets, it is proved that the model has solutions respecting
constraints on ROE's, ruin probabilities and market shares currently in practical
use. Solutions define global and optimal risk management strategies of the
company. Mathematical existence results and  tools, such as the inversion of the
linear part of the Euler-Lagrange equations, developed in a preceding paper in the
context of a simplified model are essential for the mathematical and numerical
construction of solutions of the model.  \\

\noindent {\bf Classifications JEL:} C6, G11, G22, G32 
\end{abstract}

%%%%%%%%%%%%%%%%%
\section{Introduction}
\label{Intro}
The context of this paper is a reinsurance company $H,$ with several portfolios
(or subsidiaries), being described by a Cram\`er--Lundberg like utility function
$U,$ whose value, at time $t,$ is simply the difference between the accumulated
net incomes and the accumulated claims, in the time interval $[0,t[,$
(cf. \cite{Em-Kl-Mi}).
The corporate financial problem considered is: \emph{allocate equity to the
different portfolios
(or subsidiaries), and select the portfolios, such that the annual ROE's are
satisfactory,
such that the probability of ruin of $H$ and the probability of non-solvency of
each portfolio are acceptable and such that
the expected final over-all profitability is optimal.}
We here think of a portfolio as being a portfolio of insurance
contracts or a portfolio of invested assets, which allows the inclusion of both the
liabilities and the assets in the problem.

The purpose of the present paper is three-fold: Firstly, to give a general
probabilistic set-up  of the above problem. This turns out to be possible
in terms of a stochastic optimization problem.
As we will see, mathematically, the cases with and without portfolios of
invested assets are identical, modulo a change of names of variables. We therefore
only consider the case where no invested assets are present.
The value of the utility function $U,$ at time $t,$ is then the difference between
the accumulated net premiums and the accumulated claims, in the time interval
$[0,t[.$
Secondly, to develop a method which
allows the construction of approximate solutions
satisfying the constraints of the optimization problem.
In fact, due to the non-solvency probabilities,
the problem is highly non-linear and it seems difficult to solve directly.
This is done by considering a
simplified portfolio selection model, with stronger constraints than in the
original portfolio selection model and which is easier to solve. We prove that
it is possible to choose the simplified problem so as to be the optimization of the
final expected utility under constraint on the variance of the final utility,
the expected annual ROE's and other  piece-vice linear inequality constraints. The
existence of solutions in a Hilbert space of adapted (to the claims processes)
square integrable processes, is proved under mild hypotheses on the result
processes of the portfolios.
For their construction, a Lagrangian formalism with multipliers is introduced.
Thirdly, to indicate that a basic arbitrage principle must be added to the
model, in order to eliminate a degeneracy of the allocation problem. This is
done by showing  that the initial equity allocation and the future dividends
generically are non-unique for optimal solutions, although the portfolios are
unique.

Before giving more detailed results, we shall go back to the motivation of the
problem and also introduce some notation.

The reinsurance company $H$ is organized as a holding, with subsidiaries
$S^{(1)}, \ldots, S^{(\aleph)},$ where $\aleph \geq 1$ is an integer. The companies
$S^{(i)}$ here only correspond to a division of the activities of $H$ into parts,
whose profitability, portfolio selection and certain other properties need to be
considered individually. This allows localization of
capital flows and results.
The subsidiary $S^{(i)}$ can cease its activities, which can be beneficial for $H.$
The portfolio $\theta^{(i)}$ of $S^{(i)}$ has $N^{(i)} \geq 1$ different types of
contracts $\theta^{(i)}_{j}.$ It is decomposed into
$\theta^{(i)}=\xi^{(i)}+\eta^{(i)},$ where $\xi^{(i)}$ is the run-off, at
time $t=0,$  concluded
at a finite number of past times $t < 0$ and where $\eta^{(i)}$ is the
(present and future) underwriting portfolio,
to be  concluded at a given finite number of times $0, \ldots , \bar{T},$
$\bar{T} \geq 1.$ If a subsidiary $S^{(i)}$ ceases its activities at a time
$t=t_{c},$ (past, present or future), then it only continues to
manage\label{F1}\footnote{\label{f1}The other
         possibility, which is not to keep the run-off, was considered in \cite{T1}.}
its run-off for times $t \geq t_{c}.$
The company $H$ pays dividends\footnote{Here $D(t)$ can be positive or negative.
                                   By convention, $D(t)<0$
                                        corresponds to an increase in capital.}
$D(t),$ to the shareholders at time $t.$ By convention $D(0)=0.$
Some rules determining the dividends
$D(t),$ in different situations, have been established.\label{F2}\footnote{ \label{f2} For example
       $D(t)$ can simply be a
       function of the result  $C(t)$ such that $D(t)=x(C(t)-c),$ if $C(t)
       \geq c$ and $D(t)=0,$ if $C(t) <c,$ where $0 \leq x \leq 1$ et $c \geq 0$
       are given real numbers.}
The initial equity $K(0),$ of $H$ at $t=0,$ is known.
The problem is to determine at $t=0,$
\begin{itemize}
\item i) the equity $K^{(i)}(0)$ of $S^{(i)}$ at $t=0$ ($1\leq i \leq \aleph$)
\item ii) the  dividends $D^{(i)}(n),$ which $S^{(i)}$ shall pay $H$ at time
              $n \geq 1,$
              where the determination of the dividends for time $n,$
          takes into account the observed claims during the periods preceding $n$
\item iii) the portfolio selection for $S^{(i)},$ (i.e. underwriting targets for
           present and future periods $n \geq 0$), where the determination of the
           underwriting targets $\eta^{(i)}_{j}(n)$ for the time interval $[n,n+1[,$
          takes into account the observed claims during the periods preceding $n,$
\end{itemize}
such that the constraints are satisfied and the expected utility is
optimal.

We postulate that the accessible information
$\mathcal{F}_{t}$ at time $t$ is given by the filtration generated by the
claims process.\footnote{This is an approximation, since the IBNRs at time $t$
are, by definition, not known with certainty  at $t.$}
The out-put of the model is then given by the (certain) initial equities $K^{(i)}(0),$
by the present certain $\eta^{(i)}_{j}(0)$ and future random $\eta^{(i)}_{j}(n),$
$n \geq 1$ underwriting levels
of different types of contracts and by future random dividends $D^{(i)}(n),$
$n \geq 1$ for each portfolio. By convention we set $D^{(i)}(0)=0,$ which is no
restriction.
The future random equities $K^{(i)}(n),$ $n \geq 1$
are then simply obtained by the budget constraint equations.
Points (ii) and (iii) indicate that the dividends and
the (target) underwriting
levels form stochastic processes adapted to the filtration generated by the
claims process.
The future random variables (subsequent underwriting levels and dividends) define a
strategy of reactivity to the occurrence of future exterior random events. The
uncertainty of these variables is reduced by the future increase of information. Thus
underwriting levels and dividends, at a given future time, become certain when that given
time is reached.

The probabilistic set-up of the {\em general non-linear model} is given in
\S \ref{Section 1.1}.
It extends the stochastic model first considered in \cite{T1}.
We justify, in Remark~\ref{Rm1.1.1}, that the models with and without invested
assets are mathematically identical.
The {\em general
quadratic model,} which permits a constructive approach through a Lagrangian
formalism, is given in \S \ref{Section 1.2}.
We establish (see Theorem \ref{Th1.3}), under certain mild conditions,
($h_{1}$), ($h_{2}$), ($h_{3}$) and
($h_{4}$), on the result processes for the unit-contracts,
that the optimization problem  has a solution. Theorem \ref{Th1.4} indicates
that the solution, generically, is non-unique.
Condition ($h_{1}$) says that the final
utility (sum of all results) of a unit contract, written at time $k,$
is independent of events occurring before $k$ and that the intermediate utilities
are not ``too much'' dependent. This is a starting point, since in practice, this
is generally not exactly
true, among other things because of feed-forward phenomena in the pricing.
Condition ($h_{2}$) is equivalent to the statement that no non-trivial linear
combinations of final utilities, of contracts written at a given time, is a
certain random-variable. This can also be coined, in more financial terms:
a new business portfolio (or under-writing portfolio) $\eta(t),$ constituted at
time $t,$ can not be risk-free.
Condition ($h_{3}$) says that the final
utility  of  unit contracts, written at different times are independent.
Similarly, condition ($h_{4}$) says that the final
utility  of  unit contracts, written by different subsidiaries are independent.
These conditions, which
excludes interesting situations, like cyclic markets, have only been chosen
for simplicity. They can largely be weakened without altering the results of this paper.
An important point is that no particular distributions (statistical laws) are
required.
The properties of these two models are mainly derived by considering the even
simpler model of reference \cite{T2}, here called the {\em basic model}, of which
needed facts are summed up in \S \ref{Section 1.3}.
We remind that
the portfolio in \cite{T2} is an extension of Markowitz portfolio \cite{M1}
to a multiperiod stochastic portfolio, as suggested by \cite{H-K} (c.f. also
\cite{D-J}). One of the new features is that future
results of contracts written at different times are distinguishable, which easily
allows to consider contracts with different maturity times.
The square root of the
variance of the utility of a
portfolio defines a norm, which is equivalent to the usual $L^{2}$-norm of the
portfolio (see Theorem \ref{ThA.4}).
This is one of the major technical tools in the proof of the results of the
present paper. There is existence and uniqueness
of an optimal portfolio (see Theorem \ref{ThA.4.1}).
Let us here also mention that a Lagrangian formalism is given in \cite{T2}, as well as
essential steps in the construction of the optimal solution,
(formula (2.17) of \cite{T2}).
Namely, an effective method of calculating the inverse of the linear integral operator
defined by the quadratic part of the Lagrangian, is established.
The algorithm only invokes finite dimensional linear algebra and the conditional
expectation operator.
Moreover, the determination of Lagrange multipliers is also considered in \cite{T2}.
The proofs of the results of the present article are given in \S\ref{Proof}.
For computer simulations, in the simplest cases, see \cite{Dio1}. \\

\noindent
{\bf Acknowledgement:} The author would like to thank Jean-Marie Nessi, CEO of
AXA-R\'e, and his collaborators for the many interesting discussions, which
were the starting point of this work.

\section{The model and main results}
\label{Section 1}
\subsection{General non-linear model}
\label{Section 1.1}
Mathematically, the probabilistic context of the model is given by a (separable
perfect) probability space  $(\Omega, P, \mathcal{F})$ and a filtration
$\mathcal{A}=\{\mathcal{F}_{t}\}_{t \in \mathbb{N}},$ of sub $\sigma$-algebras
of the $\sigma$-algebra $\mathcal{F},$ i.e.
$\mathcal{F}_{0}=\{\Omega,\emptyset \}$ and
$\mathcal{F}_{s} \subset \mathcal{F}_{t} \subset \mathcal{F}$ for
$0 \leq s \leq t.$ By convention $\mathcal{F}_{t}=\mathcal{F}_{0}$ for
$t <0.$

To introduce the portfolios and utility functions of the subsidiaries,
$S^{(1)},$ $\ldots,$ $S^{(\aleph)},$ let us consider the subsidiary $S^{(j)}.$
The portfolio of $S^{(j)}$ is composed of $N^{(j)} \geq 1$ types
of insurance contracts. By a unit contract, we denote a insurance contract whose
total premium is one currency unit.\footnote{All
                        flows are supposed actualized.}
The utility $u_{i}^{(j)}(t,t'),$ at $t' \in \mathbb{N}$ of the unit contract $i,$
$1 \leq i \leq N^{(j)},$ concluded at $t \in \mathbb{Z},$ is by definition
the accumulated result in the time interval $[0,t'[,$ if $t < 0 ,$
in the time interval $[t,t'[,$ if $0 \leq t \leq t'$
and $u_{i}^{(j)}(t,t')=0,$ if $0 \leq t' \leq t.$\footnote{By result we here mean the
net technical result including interest rates revenues from reserves.}
We suppose that $u_{i}^{(j)}(t,t')$
is $\mathcal{F}_{t'}$-measurable and that $(u^{(j)}(t,t'))_{t' \geq 0}$ is an element
of the space\footnote{Let $1 \leq q < \infty.$ Then 
  $(X_{i})_{0 \leq i} \in \mathcal{E}^{q}(\mathbb{R}^{N})$
  if and only if  $X_{i}:\Omega \rightarrow \mathbb{R}^{N}$ is $\mathcal{F}$
  measurable and
  $\|X_{i}\|_{L^{q}(\Omega, \mathbb{R}^{N})}
  =(\int_{\Omega} |X_{i}(\omega)|^{q}_{\mathbb{R}^{N}}dP(\omega))^{1/q}
  < \infty$
  for $i \geq 0$, where $|\thickspace \thickspace|_{\mathbb{R}^{N}}$
  is the norm in $\mathbb{R}^{N}.$  Let $\mathcal{E}^{q}
  (\mathbb{R}^{N},\mathcal{A})$ the sub-space of
   $\mathcal{A}$ adapted processes in
  $\mathcal{E}^{q}(\mathbb{R}^{N}).$ We define $\mathcal{E}(\mathbb{R}^{N})
  =\cap_{q \geq 1}\mathcal{E}^{q}(\mathbb{R}^{N})$ and
  $\mathcal{E}(\mathbb{R}^{N},\mathcal{A})
  =\cap_{q \geq 1}\mathcal{E}^{q}(\mathbb{R}^{N},\mathcal{A}).$}
$\mathcal{E}(\mathbb{R}^{N^{(j)}}),$ of processes, with finite moments of all orders.
Since, for given $t \in \mathbb{Z},$ the process $(u_{i}^{(j)}(t,t'))_{t' \geq 0}$
is $\mathcal{A}$-adapted, it follows that
$(u^{(j)}(t,t'))_{t' \geq 0} \in \mathcal{E}  (\mathbb{R}^{N^{(j)}},\mathcal{A}).$
The  final utility of the unit contract $i,$ concluded at $t,$
which is given by
$u_{i}^{(j)\infty}(t)=u_{i}^{(j)}(t,s')$ $(=u_{i}^{(j)}(t,\infty)),$ when the contract
does not generate a flow after the time $s',$ $s' \geq 0,$ is
$\mathcal{F}_{s'}$ measurable. We suppose that there exists a time $T,$
(independent of $t$)
such that the unit contracts concluded at $t \in \mathbb{Z},$ do not generate a
flow after the time $t+T.$ Let the amount of the contract
of type $i,$ where $1 \leq i \leq N^{(j)},$ concluded at time
$t \in \mathbb{Z},$
be $\theta_{i}^{(j)}(t).$ In other words,
$\theta_{i}^{(j)}(t)$ is the number of unit contracts of type $i.$ Here the run-off
$\xi^{(j)}(s)=\theta^{(j)}(s) \in \mathbb{R}^{N^{(j)}},$ $s<0$ is a certain vector
(i.e. $\mathcal{F}_{0}$-measurable) and
$\eta^{(j)}(s)=\theta^{(j)}(s),$ $s \in \{0,\ldots ,\bar{T} \}$ is a
$\mathcal{F}_{s}$-measurable random vector, taking its value in
$\mathbb{R}^{N^{(j)}}.$ We introduce a supplementary value $\tau_{f},$ which is
a final state of the process $(\theta^{(j)}(t))_{t \in \mathbb{Z}},$ reached
when the activity of the company $S^{(j)}$ ceases. In the sequel the (certain)
random variables $\xi^{(j)}(t)$ and the random variables $\eta^{(j)}(t)$ take
values in
$\mathbb{R}^{N^{(j)}} \cup \{\tau_{f}\}$\footnote{$\mathbb{R}^{N} \cup \{\tau_{f}\}$
  is not a vector space. Let
  $X_{l}: \Omega \rightarrow \mathbb{R}^{N} \cup \{\tau_{f}\},$ $l \in \{1,2\}$ 
  be two random variables and let $a \in \mathbb{R}$. We define
  $(X_{1}+X_{2})(\omega)=X_{1}(\omega)+X_{2}(\omega),$
  if $X_{l}(\omega) \neq \tau_{f}$ for $l \in \{1,2\}$ and
  $(X_{1}+X_{2})(\omega)=\tau_{f}$, if $X_{l}(\omega)= \tau_{f}$ for $l=1$ or
  $l=2.$ We also define $(aX_{1})(\omega)=a(X_{1}(\omega)),$ if
  $X_{1}(\omega) \neq \tau_{f}$ and $(aX_{1})(\omega)=\tau_{f},$ if
  $X_{1}(\omega) =\tau_{f}.$ This permits to continue
  to use the linear structure on the subspace $\mathbb{R}^{N}.$}.
Moreover it is supposed that the component of $(\eta^{(j)}(t))_{t \geq 0},$ in
$\mathbb{R}^{N^{(j)}},$ has finite variance\label{F8}\footnote{\label{f8}More precisely it is
  supposed that
  $(p \circ \eta(t))_{t \geq 0} \in \mathcal{E}^{2}(\mathbb{R}^{N},\mathcal{A}),$
  where the function $p:\mathbb{R}^{N} \cup \{\tau_{f}\} \rightarrow \mathbb{R}^{N}$
  is defined by $p(\tau_{f})=0$ and $p(x)=x,$ for $x \in \mathbb{R}^{N}.$}, i.e.
$(\eta^{(j)}(t))_{t \geq 0}
    \in \tilde{\mathcal{E}}^{2}(\mathbb{R}^{N^{(j)}} \cup \{\tau_{f}\},\mathcal{A}),$
the space of processes
with values in $\mathbb{R}^{N^{(j)}} \cup \{\tau_{f}\}$ and
and whose component in $\mathbb{R}^{N^{(j)}}$ is an element of
$\mathcal{E}^{2}(\mathbb{R}^{N^{(j)}},\mathcal{A}).$\label{F9}\footnote{\label{f9}
 A as usually a sequence
 $(X_{n})_{n \geq 1}$ in
 $\tilde{\mathcal{E}}^{2}(\mathbb{R}^{N} \cup \{\tau_{f}\},\mathcal{A})$
 converges in distributions to
 $X \in \tilde{\mathcal{E}}^{2}(\mathbb{R}^{N} \cup \{\tau_{f}\},\mathcal{A}),$
 if $E(f(X_{n}))$ converges to $E(f(X)),$ for all real bounded continuous
 functions $f$ on $\mathbb{R}^{N} \cup \{\tau_{f}\}.$
 This defines the topology of convergence in distributions (or the $d-$topology)
 on
 $\tilde{\mathcal{E}}^{2}(\mathbb{R}^{N} \cup \{\tau_{f}\},\mathcal{A}).$
 When
 a subset
 $F \subset \tilde{\mathcal{E}}^{2}(\mathbb{R}^{N} \cup \{\tau_{f}\},\mathcal{A})$
 is said to have a property of a topological vector space, such as
 being bounded, it is meant that the set
 $F' \subset \mathcal{E}^{2}(\mathbb{R}^{N},\mathcal{A}) \times
      \mathcal{E}^{2}(\mathbb{R},\mathcal{A}),$ given by
 $F'=\{(p \circ \eta, \lambda \circ \eta) \; | \; \eta \in F \},$ where
 $\lambda ( \tau_{f})=1$ and $\lambda (x)=0,$ if $x \in \mathbb{R}^{N},$ has
 that property. This convention will similarly be used for all spaces of functions
 with values in $\mathbb{R}^{N} \cup \{\tau_{f}\}.$}
We introduce the underwriting portfolio set
$\tilde{\mathcal{P}}_{u}^{(j)},$ of elements
$\eta^{(j)} \in \tilde{\mathcal{E}}^{2}(\mathbb{R}^{N^{(j)}}
                                           \cup \{\tau_{f}\},\mathcal{A}),$
such that only a finite number of $\eta^{(j)}(t),$ $t \in \mathbb{N}$ are
different from zero. Let  $\mathcal{P}_{u}^{(j)}$ be the intersection of
$\tilde{\mathcal{P}}_{u}^{(j)}$ and
$\mathcal{E}^{2}(\mathbb{R}^{N^{(j)}},\mathcal{A}).$
The spaces of run-off portfolios are defined by
$\tilde{\mathcal{P}}_{r}^{(j)}=(\mathbb{R}^{N^{(j)}} \cup \{\tau_{f}\})^{\infty}$
and $\mathcal{P}_{r}^{(j)}=(\mathbb{R}^{N^{(j)}})^{\infty}.$
If $\xi^{(j)} \in \tilde{\mathcal{P}}_{r}^{(j)}$
and if
$\eta^{(j)} \in \tilde{\mathcal{P}}_{u}^{(j)},$ then
$\theta^{(j)}=\xi^{(j)}+\eta^{(j)} \in \tilde{\mathcal{P}}^{(j)}=\tilde{\mathcal{P}}_{r}^{(j)}
 \times \tilde{\mathcal{P}}_{u}^{(j)},$ the portfolio set of $S^{(j)}.$
The aggregate portfolio $\theta =(\theta^{(1)}, \ldots, \theta^{(\aleph)}),$ of
the company $H,$ is an element of the portfolio set
$\tilde{\mathcal{P}}=\times_{1 \leq j \leq \aleph} \tilde{\mathcal{P}}^{(j)}.$
Let
$\tilde{\mathcal{P}}_{r}=\times_{1 \leq j \leq \aleph} \tilde{\mathcal{P}}_{r}^{(j)}$
and
$\tilde{\mathcal{P}}_{u}=\times_{1 \leq j \leq \aleph} \tilde{\mathcal{P}}_{u}^{(j)}.$

As already mentioned, the underwriting portfolio $\eta$ shall satisfy
constraints given by the market, the shareholders, etc. To introduce these
constraints, let $I$ be an index set and let
$G=\{g_{\alpha}| \alpha \in I \}$ be a set of functions
$g_{\alpha}:\mathbb{N} \times \tilde{\mathcal{P}} \rightarrow L^{2}(\Omega, \mathbb{R})$
such that the value $g_{\alpha}(t,\xi,\eta),$ where $\xi \in \tilde{\mathcal{P}}_{r},$
$\eta \in \tilde{\mathcal{P}}_{u},$ is independent of $\eta(t'),$ for $t' > t.$
In this paper we will say that $\eta \mapsto g_{\alpha}(t,\xi,\eta)$ is a causal
function (of $\eta$).
We define
 \begin{equation}
   \mathcal{C}(\xi,G)=\{\eta \in \tilde{\mathcal{P}}_{u}
         | g_{\alpha}(t,\xi,\eta) \geq 0, t \geq 0\},
 \label{Eq1.2}
 \end{equation}
which is the set of all underwriting portfolios $\eta$ compatibles with
the run-off $\xi \in \tilde{\mathcal{P}}_{r}$ and with the set of constraints $G.$
We note that the process $(g_{\alpha}(t,\xi,\eta))_{t \geq 0} \in 
  \mathcal{E}^{2}(\mathbb{R}).$

The utility $U^{(j)}(t, \theta^{(j)}),$ at time $t \in \mathbb{Z}$ of a portfolio
$\theta^{(j)} \in \tilde{\mathcal{P}}^{(j)},$
is defined by\footnote{Here the scalar product in $\mathbb{R}^{N}$ is
                $x \cdot y =\sum_{1 \leq i \leq N} x_{i} y_{i}.$}
\[ (U^{(j)}(t, \theta^{(j)}))(\omega)=\sum_{k \leq t}
            (\theta^{(j)}(k))(\omega) \cdot (u^{(j)}(k,t))(\omega),\]
if $(\theta^{(j)}(k))(\omega) \neq \tau_{f}$ for $k \leq t$ and by
\[ (U^{(j)}(t, \theta^{(j)}))(\omega)=\sum_{k \leq (t_{c})(\omega)-1}
            (\theta^{(j)}(k))(\omega) \cdot (u^{(j)}(k,t))(\omega),\]
if
$(t_{c})(\omega)=\inf \{k \in \mathbb{Z}  \; | \;
                  (\eta^{(j)}(k))(\omega) = \tau_{f}\} \leq t.$
The utility $U^{(j)}$ so defined can be written on the following two forms:
\begin{equation}
 \begin{split}
  U^{(j)}(t, \theta^{(j)})&
    =\sum_{k \leq t^{\ast} } \theta^{(j)}(k) \cdot u^{(j)}(k,t) \\
    &=\sum_{k \leq t}(p \circ \theta^{(j)}(k)) \cdot u^{(j)}(k,t)
              = U^{(j)}(t, p \circ \theta^{(j)}),
 \end{split}
 \label{Eq1.3}
\end{equation}
where $t \in \mathbb{Z},$ $(t^{\ast})(\omega)=\min((t_{c})(\omega)-1,t)$ and
where $p:\mathbb{R}^{N} \cup \{\tau_{f}\} \rightarrow \mathbb{R}^{N}$
is defined  by $p(\tau_{f})=0$ and $p(x)=x,$ for $x \in \mathbb{R}^{N^{(j)}},$
(c.f. footnote~\ref{f8}, p.\pageref{F8}). $U^{(j)}(t, \theta^{(j)})$ is
$\mathcal{F}_{t}$ measurable.
We have here chosen to keep the run-off for times, larger or equal to
$t_{c},$ when $S^{(j)}$ ceases its activities. Another
possibility is not to keep the run-off (c.f. footnote~\ref{f1}, p.\pageref{F1}),
in which case the utility is given by
\[\underline{U}^{(j)}(t, \theta^{(j)})=U^{(j)}(t^{\ast}, \theta^{(j)}).\]
for $t \in \mathbb{Z}.$
The  stochastic process $(U^{(j)}(t, \theta^{(j)}))_{t \geq 0}$ is an element of the
space $\mathcal{E}^{p}(\mathbb{R}, \mathcal{A}),$ for $1 \leq p < 2,$ which follows directly
from Schwarz inequality. However, without further hypotheses, it does not in
general have finite variance.
The utility of an aggregate portfolio $\theta \in \tilde{\mathcal{P}}$ is defined by
\begin{equation}
  U(t, \theta)=\sum_{1 \leq j \leq \aleph}U^{(j)}(t, \theta^{(j)}).
 \label{Eq1.3'}
 \end{equation}

The result of a portfolio 
$\theta^{(j)}
\in \tilde{\mathcal{P}}^{(j)},$
for the time period $[t,t+1[,$ $t \in \mathbb{Z},$ is defined by
\begin{equation}
  (\Delta U^{(j)})(t+1,\theta^{(j)})=\sum_{k \leq t^{\ast}} \theta^{(j)}(k) \cdot
         (u^{(j)}(k,t+1)-u^{(j)}(k,t)),
 \label{Eq1.4}
 \end{equation}
where $t^{\ast}$ is defined as in formula (\ref{Eq1.3'}).
$(\Delta U^{(j)})(t+1,\theta^{(j)})$ is $\mathcal{F}_{t+1}$ measurable. Since
$u^{(j)}(t+1,t+1)=0,$  (\ref{Eq1.3}) and (\ref{Eq1.4})
give
 \begin{equation}
  (\Delta U^{(j)})(t+1,\theta^{(j)})=U^{(j)}(t+1,\theta^{(j)})-U(t,\theta^{(j)}),
 \label{Eq1.5}
 \end{equation}
for $t \geq 0.$
The result of an aggregate portfolio $\theta \in \tilde{\mathcal{P}}$ is defined by
\begin{equation}
  (\Delta U)(t, \theta)=\sum_{1 \leq j \leq \aleph}(\Delta U^{(j)})(t, \theta^{(j)}).
 \label{Eq1.5'}
 \end{equation}

The company $S^{(j)}$ has a initial equity $K^{(j)}(0) \in \mathbb{R},$
at $t=0,$ and pays dividends $D^{(j)}(t),$ at time $t \geq 0,$ $D^{(j)}(0)=0.$ We suppose that
$(D^{(j)}(t))_{t \geq 0} \in \mathcal{E}^{2}(\mathbb{R}, \mathcal{A}).$
The dividend can be negative, which as matter of fact is an increase of equity.
The expression  (\ref{Eq1.5}), of the result for the period
$[t,t+1[,$ shows that the equity $K^{(j)}(t+1)$ at time  $t+1$ is given by
 \begin{equation}
  K^{(j)}(t+1)=K^{(j)}(t)+(\Delta U^{(j)})(t+1,\theta^{(j)})-D^{(j)}(t+1),
 \label{Eq1.6}
 \end{equation}
where $t \geq 0.$
We have that
$(K^{(j)}(t))_{t \geq 0} \in \mathcal{E}^{p}(\mathbb{R}, \mathcal{A}),$ for $1 \leq p < 2.$
The dividends  $D$ payed to the shareholders by the company $H$ and the equity
of $H$ is now given by $D=\sum_{1 \leq j \leq \aleph}D^{(j)}$ and
$K=\sum_{1 \leq j \leq \aleph}K^{(j)}$ respectively.

The companies $H$ and $S^{(1)}, \ldots, S^{(\aleph)},$ shall satisfy solvency
conditions, which are expressed as lower limits on the equity.
For a portfolio, $\theta^{(j)} \in \tilde{\mathcal{P}}^{(j)},$
let $\theta^{(j)}=\xi^{(j)}+\eta^{(j)}$ be its decomposition into
$\xi^{(j)} \in \tilde{\mathcal{P}}_{r}^{(j)}$ and $\eta^{(j)} \in \tilde{\mathcal{P}}_{u}^{(j)}.$
Moreover let
$(m^{(j)}(t,\theta^{(j)}))_{t \geq 0} \in \mathcal{E}^{2}(\mathbb{R}, \mathcal{A})$ be a
process, called solvency margin, such that
$\eta^{(j)} \mapsto m^{(j)}(t,\xi^{(j)}+\eta^{(j)})$ is a causal function and
such that $m^{(j)}(t,\theta^{(j)})=m^{(j)}(t,p \circ \theta^{(j)}),$ where
$p$ is the projection as in equation (\ref{Eq1.3}).
We define the non-solvency probability, for the portfolio
$\theta^{(j)}$ of $S^{(j)},$ with respect to the solvency margin $m^{(j)}$ by
 \begin{equation}
  \Psi^{(j)} (t,K^{(j)},\theta^{(j)},m^{(j)})
         =P(\inf\{K^{(j)}(n)-m^{(j)}(n,\theta^{(j)})|0 \leq n \leq t \}<0), 
 \label{Eq1.7}
 \end{equation}
where $t \geq 0.$ The most usual case is $m^{(j)}=0,$ i.e. positive equity,
which gives the usual ruin probability.
Similarly we define the non-solvency probability, for the portfolio
$\theta \in \tilde{\mathcal{P}}$ of $H,$ with respect to the solvency margin $m$ by
$\Psi (t,K,\theta,m)
         =P(\inf\{K(n)-m(n,\theta)|0 \leq n \leq t \}<0).$

We can now formulate the optimization problem. To precise the unknown
processes (or variables), already mentioned  in (i)-(iii) of \S\ref{Intro}, we
introduce the Hilbert space $\mathcal{P}_{u, \bar{T}}^{(j)}$
(resp. the complete metric space\label{F10}\footnote{\label{f10}The metric
     $\rho$ in $\tilde{\mathcal{P}}_{u, \bar{T}}^{(j)}$
 is defined by $\rho(X,Y)
 =(\sum_{1 \leq i \leq \bar{T}}(\|p \circ X_{i} -p \circ Y_{i} \|^{2} + P(X_{i}^{-1}(\{\tau_{f}\})\Delta Y_{i}^{-1}(\{\tau_{f}\}))^{2}))^{1/2},$
 where $\Delta$ denotes the symmetric difference of sets and $\|\thickspace \thickspace \|$
 is the norm in $L^{2}(\mathbb{R}^{N^{(j)}}).$ The topology defined by the
 metric $\rho,$ is referred to as the strong topology.}
$\tilde{\mathcal{P}}_{u, \bar{T}}^{(j)}$) of elements
$\eta^{(j)} \in \mathcal{P}_{u}^{(j)},$
(resp. $\eta^{(j)} \in \tilde{\mathcal{P}}_{u}^{(j)}$)
such that $\eta^{(j)}(t)=0,$
(resp. $(\eta^{(j)}(t))(\omega)=0$ or $(\eta^{(j)}(t))(\omega)=\tau_{f},$  $\omega \in \Omega$),
for $t > \bar{T}.$ 
We also introduce the spaces
$\mathcal{P}_{\bar{T}}^{(j)}
=\mathcal{P}_{r}^{(j)}\times \mathcal{P}_{u, \bar{T}}^{(j)},$
$\tilde{\mathcal{P}}_{\bar{T}}^{(j)}
=\tilde{\mathcal{P}}_{r}^{(j)} \times \tilde{\mathcal{P}}_{u, \bar{T}}^{(j)},$
$\mathcal{P}_{u,\bar{T}}
        =\times_{1 \leq j \leq \aleph} \mathcal{P}_{u,\bar{T}}^{(j)},$
$\tilde{\mathcal{P}}_{u,\bar{T}}
        =\times_{1 \leq j \leq \aleph} \tilde{\mathcal{P}}_{u,\bar{T}}^{(j)},$
$\mathcal{P}_{\bar{T}}=\mathcal{P}_{r} \times \mathcal{P}_{u,\bar{T}}$
and $\tilde{\mathcal{P}}_{\bar{T}}=\tilde{\mathcal{P}}_{r} \times \tilde{\mathcal{P}}_{u,\bar{T}}.$
The unknown variables are
\begin{itemize}
\item $v_{1}$)  the equity $K^{(j)}(0) \in \mathbb{R}$ of $S^{(j)}$ at $t=0$
\item $v_{2}$) the dividend process
      $D^{(j)}=(0, D^{(j)}(1), D^{(j)}(2), \ldots ) \in \mathcal{E}^{2}(\mathbb{R},\mathcal{A})$
      payed by $S^{(j)}$ to $H$
\item $v_{3}$) the underwriting portfolio $\eta^{(j)} \in \tilde{\mathcal{P}}_{u, \bar{T}}^{(j)},$
       of $S^{(j)},$ with the underwriting horizon $\bar{T} \in \mathbb{N}+1$
       fixed (independent of $j$),
\end{itemize}
where $1 \leq j \leq \aleph.$
 Thus the unknown variables of the subsidiary $S^{(j)}$
are the components of the variable
 $Z^{(j)}=(K^{(j)}(0),D^{(j)}, \eta^{(j)}(0), \ldots ,\eta^{(i)}(\bar{T}),
\eta^{(j)}(\bar{T}+1),\ldots),$ where $\eta^{(j)}(\bar{T}+l)$ can only take
the values $0$ and $\tau_{f},$ for $l > 0.$
The optimization is done with respect to the variable
$Z=(Z^{(1)}, \ldots ,Z^{(\aleph)}),$ satisfying ($v_{1}$)-($v_{3}$).
It will be convenient to use a variable, which is obtained by a permutation of
the coordinates of $Z.$ Let $\Vec{K}(0)=(K^{(1)}(0), \ldots ,K^{(\aleph)}(0))$
and $\Vec{D}=(D^{(1)}, \ldots ,D^{(\aleph)}).$ We use the variable
$(\eta, \Vec{K}(0), \Vec{D})$ instead of $Z.$

Certain data are given:
\begin{itemize}
\item $d_{1}$) initial equity $K(0) \in \mathbb{R}^{+}$ of $H$
\item $d_{2}$) the dividend process $D(\theta) \in \mathcal{E}^{2}(\mathbb{R},\mathcal{A}),$
				with $(D(\theta))(0)=0,$ which $H$ pays the share holders, only
				depends of, past and present, aggregate results, (i.e.
				$(D(\theta))(t)=f_{t}((\Delta U)(1, \theta), \ldots, (\Delta U)(t, \theta)),$
				$t \geq 1,$ where $f_{t}$ is a $\mathcal{F}_{t}$-measurable function,
				c.f. footnote~\ref{f2}, p.\pageref{F2})
\item $d_{3}$) the run-off $\xi \in \tilde{\mathcal{P}}_{r}$
               of the subsidiary.
\end{itemize}

Before introducing the constraints, we recall that no flows are generate by
contracts after a certain time $\bar{T}+T.$ Therefore equity $K^{(j)}(t)$ is
constant for $t \geq \bar{T}+T.$ We suppose that the solvency margins are also
chosen such that they are constant for sufficiently large times. We choose $T$
such that they are constant for $t \geq \bar{T}+T.$ The constraints are:
\begin{itemize}
\item $c_{1}$) $K(0)=\sum_{1 \leq j \leq \aleph}K^{(j)}(0)$
\item $c_{2}$) $D=\sum_{1 \leq j \leq \aleph}D^{(j)}$
\item $c_{3}$) $E((\Delta U)(t+1, \xi+\eta))
				\geq c(t)E(\sum_{1 \leq j \leq \aleph}K^{(j)}(t)),$ where
				$c(t) \in \mathbb{R}^{+}$ is given for $t \in \mathbb{N},$ 
				(constraint on ROE)
\item $c_{4}$) $\Psi(t,\sum_{1 \leq j \leq \aleph}K^{(j)},\xi+\eta,0)
				\leq \epsilon (t),$ where
				$\epsilon(t) \in [0,1]$ is a given acceptable ruin probability of
				$H,$ for $t \in \mathbb{N}$
\item $c_{5}$) supplementary constraints, to be specified, on $D^{(j)},$
				(e.g. to increase the equity of $S^{(l)},$ one can set $D^{(l)}=0$).
				To be general, we only suppose that there are real valued functions
				$F_{\alpha},$ $\alpha \in I,$ an index set, such that
				$F_{\alpha}(t,\eta, \Vec{K}(0), \Vec{D}) \leq C_{\alpha}(t,K(0),\xi),$
				where $C_{\alpha}(t,K(0),\xi)$ are constants only depending of
				the initial equity and the run-off.
\item $c_{6}$)  if $(\eta^{(j)}_{i}(t))(\omega)
				\neq \tau_{f},$ then $((\underline{c}^{(j)}_{i}(\eta))(t))(\omega)  \leq
				(\eta^{(j)}_{i}(t))(\omega) < \infty$ and
				$(\eta^{(j)}_{i}(t))(\omega) \leq ((\overline{c}^{(j)}_{i}(\eta)))(t))(\omega),$
				for $\omega \in \Omega$ (a.e.), $1 \leq j \leq \aleph,$
				$1 \leq i \leq N^{(j)}$ and $0 \leq t \leq \bar{T}.$
				Here
				$\underline{c}^{(j)}_{i}(\eta) \in \mathcal{E}^{2}(\mathbb{R},\mathcal{A})$
				and  $\overline{c}^{(j)}_{i}(\eta)$
				are given $\mathcal{A}$-adapted processes, which are causal functions
				of $\eta$ and which satisfy $((\underline{c}^{(j)}_{i}(\eta))(t))(\omega) \in [0,\infty[$
				and $((\overline{c}^{(j)}_{i}(\eta))(t))(\omega) \in [0,\infty],$
				(market constraints)
\item $c_{7}$) $\Psi^{(j)} (t,K^{(j)},\xi^{(j)}+\eta^{(j)},
				m^{(j)}) \leq \epsilon^{(j)}(t),$ where
				$(m^{(j)}(t,\xi^{(j)}+\eta^{(j)}))_{t \geq 0}
				\in \mathcal{E}^{2}(\mathbb{R},\mathcal{A})$ is a given solvency margin and
				$\epsilon^{(j)}(t) \in [0,1]$ is a given acceptable non-solvency
				probability of $S^{(j)},$ for $1 \leq j \leq \aleph$ and $t \in \mathbb{N}$
\item $c_{8}$) if $t > (t^{(j)}_{f})(\omega),$
				then $(\eta^{(j)}(t))(\omega)=\tau_{f},$ where $(t^{(j)}_{f})(\omega)$
				is the smallest time in $\mathbb{N}$ such that $K^{(j)}(t^{(j)}_{f})
				< m^{(j)}(t^{(j)}_{f},\xi^{(j)}+\eta^{(j)}),$
				and if $t \leq (t^{(j)}_{f})(\omega)$ and
				$(K^{(j)}(t-1))(\omega)> (m^{(j)}(t-1,\xi^{(j)}+\eta^{(j)}))(\omega),$
				then $(\eta^{(j)}(t))(\omega) \in \mathbb{R}^{N^{(j)}},$ for
				$\omega \in \Omega$ (a.e.), $1 \leq j \leq \aleph$ and $t\geq 1$
				(the activity of $S^{(j)}$ ceases just after that
				the solvency margin is not satisfied)
\end{itemize}
Let $\mathcal{C}_{c}$ be the set of all $(\eta, \Vec{K}(0), \Vec{D}),$
satisfying ($v_{1}$)-($v_{3}$) and satisfying
the constraints ($c_{1}$)--($c_{8}$).
Thus we sum up the constraints ($c_{1}$)--($c_{8}$) on the form:
\begin{itemize}
\item $c$) $(\eta, \Vec{K}(0), \Vec{D}) \in \mathcal{C}_{c}$
\end{itemize}

Among all the possible functions to optimize, we simply choose the expected
value of the final utility, $\eta \mapsto E(U(\infty, \eta+\xi)).$ The
optimization problem is now: given the initial equity $K(0),$
the dividend process $D(\xi+\eta),$ as a function of $\eta,$ and the run-off $\xi$
of $H,$ satisfying ($d_{1}$,) ($d_{2}$) and ($d_{3}$) respectively, find the solutions
$(\widehat{\eta}, \widehat{\Vec{K}}(0), \widehat{\Vec{D}}) \in \mathcal{C}_{c}$
of the equation
 \begin{equation}
   E(U(\infty, \widehat{\eta}+\xi))=\sup_{(\eta, \Vec{K}(0), \Vec{D}) \in
         \mathcal{C}_{c}} E(U(\infty, \eta+\xi)).
 \label{Eq1.8}
 \end{equation}
Due to the constraints ($c_{4}$) on the ruin probability and ($c_{7}$) on the
non-solvency probabilities, the resolution of this optimization problem leads
to highly non-linear equations. This is true even in the case of practical
applications, where the other constraints usually are  piece-vice linear.

\begin{remark} \label{Rm1.1}
The constraints ($c_{1}$) and ($c_{2}$) are just budget constraints. We have
chosen the simplest form of the constraints ($c_{3}$) on the ROE. Another
possibility, is to strengthen it so that the ROE
for the time interval $[t,t+1[$ is satisfied conditionally to the information
available at time $t,$ i.e. $E((\Delta U)(t+1, \xi+\eta)| \mathcal{F}_{t}) \geq
c(t)\sum_{1 \leq j \leq \aleph}K^{(j)}(t).$ Of course, the
expected value of the final utility will then in general be smaller for a
optimal solution.  ($c_{4}$) is just  a ruin constrain for $H.$ ($c_{5}$) is a
very general constraint, which should cover most cases coming up in applications.
($c_{6}$)
says that, if $S^{(j)}$ is in business at time $t,$ then the underwriting targets
for t, $(\eta^{(j)}_{i}(t))(\omega)$ is in a given semi bounded or bonded closed
interval of positive numbers. Often the limits are proportional to
$(\eta^{(j)}_{i}(t-1))(\omega).$ In constraint ($c_{8}$), $(t^{(j)}_{f})(\omega)$
is the smallest time such that the solvency margin is strictly negative for
$S^{(j)}$ in the state $\omega.$ The constraint says that $S^{(j)}$ ceases its
activities for times larger than $(t^{(j)}_{f})(\omega).$ More-over it says that,
if the solvency margin is strictly positive for a time $t$ before
$(t^{(j)}_{f})(\omega),$ then $S^{(j)}$ does not cease its activities at $t+1.$
So, it is only when the solvency margin is zero, there is a choice.
\end{remark}

The constraints ($c_{1}$)--($c_{5}$) and ($c_{7}$) have a form as in
formula (\ref{Eq1.2}).
The constraints ($c_{6}$) and ($c_{8}$) can also be written on this form, which
we give for later reference.
Let $\lambda_{j} : \mathbb{R}^{N^{(j)}} \cup \{\tau_{f}\} \mapsto \{0,1\}$
be defined by $\lambda_{j}(\tau_{f})=1$ and $\lambda_{j}(x)=0,$
$x \in \mathbb{R}^{N^{(j)}}.$ The constraint ($c_{6}$) is then given by
\begin{equation}
(1-\lambda_{j} \circ \eta^{(j)}(t))(\eta^{(j)}_{i}(t)-(\underline{c}^{(j)}_{i}(\eta))(t)) \geq 0
 \label{Eq1.8.0.2}
\end{equation}
and
\begin{equation}
(1-\lambda_{j} \circ \eta^{(j)}(t))((\overline{c}^{(j)}_{i}(\eta))(t))-\eta^{(j)}_{i}(t)) \geq 0,
 \label{Eq1.8.0.3}
\end{equation}
where $1 \leq j \leq \aleph,$ $1 \leq i \leq N^{(j)}$ and $0 \leq t \leq \bar{T}+T.$
In the case of ($c_{8}$) let
$s^{(j)}(t,x)=\min_{0 \leq k \leq t} (K^{(j)}(k)- m^{(j)}(k,\xi^{(j)}+\eta^{(j)})),$
where $K^{(j)}(k)$ is evaluated at $x =(\eta, \Vec{K}(0), \Vec{D}).$ Let the
step function $H$ be defined by $H(s)=0$ if $s<0$ and $H(s)=1$ if $0 \leq s.$ The
constraint ($c_{8}$) is then given by
\begin{equation}
(1-\lambda_{j} \circ \eta^{(j)}(t))s^{(j)}(t-1,x)\geq 0,
 \label{Eq1.8.0.4}
\end{equation}
\begin{equation}
\lambda_{j} \circ \eta^{(j)}(t)s^{(j)}(t-1,x)\leq 0
 \label{Eq1.8.0.5}
\end{equation}
and
\begin{equation}
H(s^{(j)}(t-1,x))(\lambda_{j} \circ \eta^{(j)}(t))(K^{(j)}(t-1)-m^{(j)}(k,\xi^{(j)}+\eta^{(j)}))= 0,
 \label{Eq1.8.0.6}
\end{equation}
where $1 \leq j \leq \aleph,$ $1 \leq i \leq N^{(j)}$ and $1 \leq t \leq \bar{T}+T+1.$

\begin{remark} \label{Rm1.1.1}
We shall illustrate that the mathematical formalism, in the cases with and without
invested assets, are identical. Let $\eta=(\eta',\eta^{(\aleph)})$ and let
$u^{(\aleph)}(t,t')=p(t+1)-p(t),$ for $t<t',$ where
$p \in \mathcal{E}(\mathbb{R}^{N^{(\aleph)}},\mathcal{A}),$ $p(0)=0$ and
$p(t)=p(\bar{T}),$ for $t \geq \bar{T}.$
Then
\begin{equation}
   U(t, \eta)=U'(t, \eta') +\sum_{0 \leq k \leq t-1} \eta^{(\aleph)}(k) \cdot (p(k+1)-p(k)),
 \label{Eq1.8.1}
\end{equation}
where $U'$ is the utility function given by $U'(t, \eta')=U(t, (\eta',0)),$
i.e. for the portfolios $1,2, \ldots, \aleph-1.$ The second term on the right
hand side of (\ref{Eq1.8.1}) is exactly the accumulated income, in the time interval
$[0,t[,$ from a portfolio $\eta^{(\aleph)}$ of invested assets, with price $p(t)$
at time $t.$ This shows that the invested assets are taken into account by the
model of this paper, as a particular case.
\end{remark}

\subsection{General quadratic model}
\label{Section 1.2}
The constraints ($c_{4}$)
on the ruin probability and ($c_{7}$) on the non-solvency probabilities are
replaced by \textit{stronger} quadratic constraints, in this model. It will be
supposed that the non-ruin and non-solvency margins are satisfied in the mean.
We introduce the constraints, where $V$ denote the variance operator:
\begin{itemize}
\item $c_{4}'$) $V(\sum_{1 \leq j \leq \aleph}K^{(j)}(t))
				\leq \epsilon' (t) (\delta(t)K(0))^{2}$ and
				$E(\sum_{1 \leq j \leq \aleph}K^{(j)}(t)) \geq \delta(t)K(0),$
				$t \in \mathbb{N},$
				where $\epsilon' (t) \geq 0$ and  $\delta(t) > 0$
\item $c_{7}'$) $V(K^{(j)}(t)-m^{(j)}(t,\xi^{(j)}+\eta^{(j)}))
				\leq {\epsilon'}^{(j)}(t) (\delta^{(j)}(t)K(0))^{2}$
				and
				$E(K^{(j)}(t)$ $-m^{(j)}(t,\xi^{(j)}+\eta^{(j)})) \geq \delta^{(j)}(t)K(0),$
				for $1 \leq j \leq \aleph$ and $t \in \mathbb{N},$
				where $m^{(j)}$ are as in ($c_{7}$)  and where
				${\epsilon'}^{(j)}(t) \geq 0$ and $\delta^{(j)}(t) >0$
\end{itemize}
For given initial equity $K(0),$ dividend process $D(\xi+\eta),$ as a function of
$\eta,$ and run-off $\xi$
of $H,$ satisfying ($d_{1}$), ($d_{2}$) and ($d_{3}$) respectively,
let $\mathcal{C}_{c'}$ be the set of variables $(\eta, \Vec{K}(0), \Vec{D}),$
satisfying ($v_{1}$), ($v_{2}$) and ($v_{3}$)  and satisfying  the constraints
($c_{1}$)--($c_{3}$), ($c'_{4}$), ($c_{5}$), ($c_{6}$), ($c'_{7}$) and ($c_{8}$).
We sum up the constraints of the quadratic model on the form:
\begin{itemize}
\item $c'$) $(\eta, \Vec{K}(0), \Vec{D}) \in \mathcal{C}_{c'}$
\end{itemize}
\begin{remark} \label{Rm1.2}
Of-course this model is not quadratic in $\eta$ for general $m,$ $F_{\alpha},$
$\underline{c}$ and $\overline{c}.$ However for common choices of these functions
it is piece-vice quadratic, which is the reason to keep the name quadratic.
\end{remark}
The optimization problem, in the case of the quadratic model, is now:
given the initial equity $K(0),$
the dividend process $D(\xi+\eta),$ as a function of $\eta,$ and the run-off $\xi$
of $H,$ satisfying ($d_{1}$), ($d_{2}$) and ($d_{3}$) respectively, find the solutions
$(\widehat{\eta}, \widehat{\Vec{K}}(0), \widehat{\Vec{D}}) \in \mathcal{C}_{c'}$
of the equation
 \begin{equation}
   E(U(\infty, \widehat{\eta}+\xi))=\sup_{(\eta, \Vec{K}(0), \Vec{D}) \in
         \mathcal{C}_{c'}} E(U(\infty, \eta+\xi)).
 \label{Eq1.9}
 \end{equation}
The constraints in the quadratic optimization problem (\ref{Eq1.9}) are stronger
than those in the original problem (\ref{Eq1.8}).
\begin{theorem} \label{Th1.2}
 If $\sum_{0 \leq k \leq t}\epsilon' (k) \leq \epsilon (t)$ and
 $\sum_{0 \leq k \leq t}{\epsilon'}^{(j)}(k) \leq {\epsilon}^{(j)}(t),$
 for $t \in \mathbb{N}$ and $1 \leq j \leq \aleph,$  then
 $\mathcal{C}_{c'} \subset \mathcal{C}_{c}.$
\end{theorem}
In order to give, in this paper, a mathematical analysis, which is as simple as
possible, of optimization problem (\ref{Eq1.9}), we shall make certain (technical)
hypotheses on the claim processes. The following hypotheses give a clear-cut
mathematical context:
\begin{itemize}
\item $h_{1}$) independence with respect to the past:
	\begin{enumerate}
			\item $u^{(p)\infty}(k)$ is independent of $\mathcal{F}_{k}$ for
				$k \in \mathbb{Z}$ and $1 \leq p \leq \aleph$
			\item  $\|E((u^{(p)}(k,t))^{2}| \mathcal{F}_{k}) \|_{L^{\infty}} < \infty$
				for $k < t$
	\end{enumerate}
\item $h_{2}$) for $k \in \mathbb{Z}$ and $1 \leq p \leq \aleph,$
				the  $N  \times N$ (positive) matrix
				$c^{(p)}(k)$ with elements
				$c^{(p)}_{ij}(k)=E((u_{i}^{(p)\infty}(k)-E(u_{i}^{(p)\infty}(k))
				(u_{j}^{(p)\infty}(k)-E(u_{j}^{(p)\infty}(k)))$ is strictly positive
\item $h_{3}$) $u_{i}^{(p)\infty}(k)$ and $u_{j}^{(p)\infty}(l)$ are independent for
				$k \neq l$
\item $h_{4}$) $u_{i}^{(p)\infty}(k)$ and $u_{j}^{(r)\infty}(l)$ are independent for
				$p \neq r$
\end{itemize}
We note that the second point of ($h_{1}$) is trivially satisfied if
$u^{(p)}(k,t)$ is independent of $\mathcal{F}_{k},$ for $k < t.$

Next theorem give the existence of optimal solutions of problem (\ref{Eq1.9}).
Approximations of these solutions can be constructed, using a Lagrangian
formalism, as in the case of the basic model in \S \ref{Section 1.3}. In order
to state the theorem, we remind that if $\xi$ is as in (\ref{Eq1.9}), then the
functions $\eta^{(j)} \mapsto m^{(j)}(t,\xi^{(j)}+\eta^{(j)}),$
$\eta \mapsto D(\xi +\eta),$ $\eta \mapsto \underline{c}^{(j)}_{i}(\eta),$
$\eta \mapsto \overline{c}^{(j)}_{i}(\eta)$ and
$(\eta, \Vec{K}(0), \Vec{D}) \mapsto
      F_{\alpha}(t,\eta, \Vec{K}(0), \Vec{D})$
are defined for
$\eta^{(j)} \in \tilde{\mathcal{P}}_{u, \bar{T}}^{(j)},$
$\eta \in \tilde{\mathcal{P}}_{u, \bar{T}}$ and
$(\eta, \Vec{K}(0), \Vec{D}) \in \tilde{\mathcal{P}}_{u, \bar{T}}
 \times {\mathbb{R}}^{\aleph}
 \times \mathcal{E}^{2}({\mathbb{R}}^{\aleph},\mathcal{A}).$
We also remind that formulas (\ref{Eq1.5}), (\ref{Eq1.5'}) and (\ref{Eq1.6})
give
 \begin{equation}
  K^{(j)}(t)=K^{(j)}(0)+ U^{(j)}(t,\theta^{(j)})-\sum_{1 \leq k \leq t}D^{(j)}(k),
 \label{Eq1.10}
 \end{equation}
$\theta^{(j)} \in \tilde{\mathcal{P}}_{ \bar{T}}^{(j)},$
and that formula (\ref{Eq1.3'}) gives
 \begin{equation}
  K(t)=K(0)+ U(t,\theta)-\sum_{1 \leq k \leq t}D(k,\theta),
 \label{Eq1.11}
 \end{equation}
when constraints ($c_{1}$) and ($c_{2}$) are satisfied and
$\theta \in \tilde{\mathcal{P}}_{ \bar{T}}.$
\begin{theorem} \label{Th1.3}
 Let the utilities $u^{(p)}(k,t)$ and $u^{(p)\infty}(k),$ of unit contracts, satisfy
 ($h_{1}$)--($h_{4}$). Let the functions
 $\eta^{(j)} \mapsto m^{(j)}(t,\xi^{(j)}+\eta^{(j)}),$
 $\eta \mapsto D(\xi +\eta)$ and $\eta \mapsto \underline{c}^{(j)}_{i}(\eta)$
 to $\mathcal{E}^{2}(\mathbb{R},\mathcal{A})$
 map bounded sets into bounded sets. In the d-topology, let
 $\eta^{(j)} \mapsto m^{(j)}(t,\xi^{(j)}+\eta^{(j)})$ and $\eta \mapsto D(\xi +\eta)$
 be continuous,  let
 $(\eta, \Vec{K}(0), \Vec{D}) \mapsto F_{\alpha}(t,\eta, \Vec{K}(0), \Vec{D}),$
 $\alpha \in I,$ be lower semi continuous and let
 $\eta \mapsto (\underline{c}^{(j)}_{i}(\eta))(t,\omega)$ and
 $\eta \mapsto -(\overline{c}^{(j)}_{i}(\eta))(t,\omega)$ be lower semi continuous (a.e.).
 If \[V(\sum_{1 \leq k \leq \bar{T}+T} D(k,\theta)) \leq c^{2} V(U(\infty, \theta)),\] where
 $0 \leq c < 1,$ and if $\mathcal{C}_{c'}$ is non-empty, then the optimization
 problem (\ref{Eq1.9}) has a solution
 $\widehat{x}=(\widehat{\eta}, \widehat{\Vec{K}}(0), \widehat{\Vec{D}}) \in \mathcal{C}_{c'}.$
\end{theorem}
\begin{remark} \label{Rm1.3}
In applications it is easy to verify that the boundedness and continuity properties
are satisfied. The variance condition simply translates that the final accumulated
dividend are less volatile than the the final accumulated result.
\end{remark}
\begin{remark} \label{Rm1.3.1}
The optimization problem (\ref{Eq1.9}) can be formulated as the optimization
of a Lagrangian with multipliers. To show this let $\lambda_{j}$ be as
in (\ref{Eq1.8.0.2}), let $p_{j}:\mathbb{R}^{N^{(j)}} \cup \{\tau_{f}\} \rightarrow \mathbb{R}^{N^{(j)}}$
be defined by $p_{j}(\tau_{f})=0$ and $p_{j}(x)=x,$ for $x \in \mathbb{R}^{N^{(j)}},$
let $p=(p_{1}, \ldots ,p_{\aleph})$ and $\lambda=(\lambda_{1}, \ldots ,\lambda_{\aleph}).$
Let $\mathcal{E}^{2}_{\bar{T}}(\mathbb{R}^{N},\mathcal{A})$ be the subset of
elements $\eta \in \mathcal{E}^{2}(\mathbb{R}^{N},\mathcal{A}),$ such that
$\eta(t) =0$ for $t>\bar{T}.$ By definition, if $\eta \in \tilde{\mathcal{P}}_{u,\bar{T}}$
then $p \circ \eta \in \mathcal{E}^{2}_{\bar{T}}(\mathbb{R}^{N},\mathcal{A}),$ where
$N=\sum_{1 \leq j \leq \aleph} N^{(j)}$ and
$\lambda \circ \eta \in \mathcal{E}^{2}_{\bar{T}}(\mathbb{R}^{\aleph},\mathcal{A}).$ We
introduce the Hilbert space
$H_{0}=\mathcal{E}^{2}_{\bar{T}}(\mathbb{R}^{N},\mathcal{A})
 \oplus \mathcal{E}_{\bar{T}}^{2}(\mathbb{R}^{\aleph},\mathcal{A})$ and the
Hilbert space $H_{1}$ of elements
$( \Vec{K}(0), \Vec{D}) \in \mathbb{R}^{\aleph}
            \oplus \mathcal{E}_{\bar{T}}^{2}(\mathbb{R}^{\aleph},\mathcal{A}),$
such that $\Vec{D}(0)=0.$ Let $H=H_{0} \oplus H_{1}.$
The optimization problem can now formulated in the variable
$x=(\alpha, \beta, \Vec{K}(0), \Vec{D})$ in $H,$ where
$(\alpha, \beta) \in H_{0}$ and where for solutions $\alpha=p \circ \eta$ and
$\beta=\lambda \circ \eta.$ The constraints $(1-\beta^{(j)})\beta^{(j)}=0$
and $\alpha^{(j)}\beta^{(j)}=0$ shall then be satisfied. The constraints
($c_{1}$)--($c_{5}$), ($c_{4}'$) and ($c_{7}'$) are easily expressed in the new
variables. The constraint ($c_{6}$) is reformulated by using (\ref{Eq1.8.0.2})
and (\ref{Eq1.8.0.3}) and constraint ($c_{8}$) by (\ref{Eq1.8.0.4})--(\ref{Eq1.8.0.6}).
This gives a Lagrangian with multipliers. We note that the function on the left hand side
of (\ref{Eq1.8.0.6}) is not differentiable, which leads to singularities in
the Euler-Lagrange equation. Approximation schemes can be based on the inversion
methods developed in \cite{T2}, for the linear part of the Euler-Lagrange equation.
A detailed solution of this problem is the subject of future studies.
\end{remark}
In general the solution $\widehat{x} \in \mathcal{C}_{c'},$ of Theorem \ref{Th1.3},
is not unique. This fact can be traced back to a simplified case, namely where
only the constraints ($c_{1}$)--($c_{4}$) are considered and where all
$\eta^{(j)}_{i} \geq 0.$ To state the result let $\mathcal{C}_{c''}$ be the set
of all $(\eta, \Vec{K}(0), \Vec{D}) \in \tilde{\mathcal{P}}_{u, \bar{T}}
 \times {\mathbb{R}}^{\aleph}
 \times \mathcal{E}^{2}({\mathbb{R}}^{\aleph},\mathcal{A}),$ satisfying
($v_{1}$)--($v_{3}$), ($c_{1}$), ($c_{2}$) and ($c'_{4}$).
We consider the following optimization problem: given $K(0) \geq 0,$ $D=0$ and
$\xi=0,$ find the solutions $\widehat{x} \in \mathcal{C}_{c''}$ of the equation
 \begin{equation}
  E(U(\infty, \widehat{\eta}))=\sup_{(\eta, \Vec{K}(0), \Vec{D}) \in
         \mathcal{C}_{c''}} E(U(\infty, \eta)).
 \label{Eq1.12}
 \end{equation}
\begin{theorem} \label{Th1.4}
If $\mathcal{C}_{c''}$ is non-empty, then the optimization problem (\ref{Eq1.12})
has a solution $\widehat{x} \in \mathcal{C}_{c''}.$ Moreover, $\widehat{\eta}$
is unique and $(\widehat{\Vec{K}}(0), \widehat{\Vec{D}})$ only has to satisfy
$K(0)=\sum_{1 \leq j \leq \aleph}  \widehat{K^{(j)}}(0)$ and
$0=\sum_{1 \leq j \leq \aleph}  \widehat{D^{(j)}}.$
\end{theorem}
In the situation of the theorem, there is a whole hyperplane
( ${\mathbb{R}}^{\aleph-1}
 \times \mathcal{E}^{2}({\mathbb{R}}^{\aleph-1},\mathcal{A})$ translated)
of solutions in the variable $(\widehat{\Vec{K}}(0), \widehat{\Vec{D}}).$
The generic case seems to be close to this case. To avoid a heavy ``book-keeping''
of solutions we only illustrate this by an informal remark instead of stating a
formal theorem.
\begin{remark} \label{Rm1.4}
Suppose that a solution $\widehat{x}$ of Theorem \ref{Th1.4} satisfies condition
($c_{7}'$) for given $m^{(j)}$ and $\epsilon^{(j)},$ $1 \leq j \leq \aleph.$
Apart from exceptional cases, there will be a whole neighbourhood of elements
$(\Vec{K}'(0), \Vec{D}')$ in a submanifold homeomorphic to ${\mathbb{R}}^{\aleph-1}
 \times \mathcal{E}^{2}({\mathbb{R}}^{\aleph-1},\mathcal{A}),$
which also give solutions $y=(\widehat{\eta},\Vec{K}'(0), \Vec{D}').$ In the
general situation of Theorem \ref{Th1.2}, this degeneracy can be partially reduced
because some of the constraints will saturate. However a supplementary economic
principle seems to be needed in order to guarantee uniqueness. If not guaranteed,
then it does not always matter what we do with the equity!
\end{remark}

\subsection{Basic model}
\label{Section 1.3}
We shall here sum up certain results obtained in reference \cite{T2} concerning
a particularly simple model, which is an essential building
block of the already considered general models. In that model it is supposed that the number
of subsidiaries $\aleph =1,$ the run-off $\xi=0,$ the dividends $D=0$  and it is
supposed that there are
no market limitations on the subscription levels, except that they are positive.
It is also imposed that the portfolio $\eta$ is an element of the Hilbert space
$\mathcal{H}=\mathcal{P}_{u, \bar{T}}.$ (So
$\eta(\omega) \neq \tau_{f}$ on $\Omega$).
We remind that, in this situation, the equity
\begin{equation}
K(t)=K(0)+U(t,\eta),
 \label{EqA.-2N}
 \end{equation}
where $K(0) \geq 0$ is the initial equity at $t=0.$

In the sequel of this paragraph, we closely follow reference \cite{T2}.
Constraints on the variable $\eta$ are introduced:
\begin{itemize}
\item $C_{3}$) $E((\Delta U)(t+1,\eta))
				\geq c(t)E(K(t)),$
				$c(t) \in {\mathbb R}^{+}$ is given (constraint on profitability)
\item $C_{4}$) $E(( U(\infty, \eta)
				-E(U(\infty, \eta)))^{2})
				\leq \sigma^{2},$ where $\sigma^{2} > 0$ is given
				( acceptable level of the variance of the final utility)
\item $C_{6}$) $0 \leq \eta_{i}(t),$ where $1 \leq i \leq N$
				(only positive subscription levels)
\end{itemize}
Let $\mathcal{C}_{0}$ be the set of portfolios $\eta \in \mathcal{H}$ such that
constraints ($C_{3}$), ($C_{4}$) and ($C_{6}$) are satisfied. This is well-defined.
In fact the quadratic form
\begin{equation}
 \eta \mapsto \mathfrak{a}(\eta)= E(( U(\infty,\eta))^{2}),
 \label{EqA.-1N}
 \end{equation}
in $\mathcal{H},$ has a maximal domain $\mathcal{D}(\mathfrak{a}),$
since for each $\eta \in \mathcal{H},$
the  stochastic process $(U(t, \theta)_{t \geq 0}$ is an element of the
space $\mathcal{E}^{p}(\mathbb{R}, \mathcal{A}),$ for $1 \leq p < 2$ (which
follows directly from Schwarz inequality).
The optimization
problem  is now, to find all
$\hat{\eta} \in \mathcal{C}_{0},$ such that
\begin{equation}
  E(U(\infty, \hat{\eta}))
   =\sup_{\eta \in \mathcal{C}_{0}} E(U(\infty, \eta)).
 \label{EqA.0}
\end{equation}
The solution of this optimization problem is largely based on the study of the
quadratic form
\begin{equation}
 \eta \mapsto \mathfrak{b}(\eta)= E(( U(\infty,\eta) -E(U(\infty,\eta)))^{2}),
 \label{EqA.1N}
 \end{equation}
in $\mathcal{H},$ with (maximal) domain
$\mathcal{D}(\mathfrak{b})=\mathcal{D}(\mathfrak{a}).$

We make certain (technical) hypotheses on the claim processes:
\begin{itemize}
\item $H_{1}$) $u^{\infty}(k)$ is independent of $\mathcal{F}_{k}$ for
				$k \in \mathbb{N}$
\item $H_{2}$) for $k \in \mathbb{N}$ the  $N  \times N$ (positive) matrix
				$c(k)$ with elements
				$c_{ij}(k)=E((u_{i}^{\infty}(k)-E(u_{i}^{\infty}(k))
				(u_{j}^{\infty}(k)-E(u_{j}^{\infty}(k)))$ is strictly positive
\item $H_{3}$) $u_{i}^{\infty}(k)$ and $u_{j}^{\infty}(l)$ are independent for
				$k \neq l$
\end{itemize}
%%%%%%%%%%%%%%%%%%
The next crucial result (Lemma 2.2 and Theorem 2.3 of \cite{T2})
show that (the square root of) each one of the quadratic forms
$\mathfrak{b}$ and $\mathfrak{a}$ is equivalent to the norm in $\mathcal{H}.$
%%%%%%%%%%%%%%%%%%%%%%%%%%%%%%%%%%%%%%%%%%%
\begin{theorem}\label{ThA.4}
If the hypotheses ($H_{1}$), ($H_{2}$) and ($H_{3}$) are satisfied, then the
quadratic forms $\mathfrak{b}$ and $\mathfrak{a}$
are bounded from below and from above,  by  strictly positive numbers $c$ and $C$ respectively,
where $0 < c \leq C,$  i.e.
\begin{equation}
   c \| \eta \|^{2}_{\mathcal{H}} \leq \mathfrak{b}(\eta) \leq
       \mathfrak{a}(\eta) \leq C \| \eta \|^{2}_{\mathcal{H}},
 \label{EqA.4.1}
 \end{equation}
for $\eta \in \mathcal{H}.$
\end{theorem}
%%%%%%%%%%%%%%%%%%%%%%%
The operators $B$ (resp. $A$) in
$\mathcal{H},$ associated with $\mathfrak{b}$ (resp. $\mathfrak{a}$),
(by the representation theorem),
i.e.
\begin{equation}
  \mathfrak{b}(\xi,\eta)=(\xi,B \eta)_{\mathcal{H}}
 \quad (\text{resp.} \quad
  \mathfrak{a}(\xi,\eta)=(\xi,A \eta)_{\mathcal{H}}),
 \label{EqA.5}
 \end{equation}
for $\xi \in \mathcal{H}$ and $\eta \in \mathcal{H},$
are strictly positive, bounded, self-adjoint operators onto ${\mathcal{H}}$
with bounded inverses.
There exist $c \in \mathbb{R},$
such that $0 < cI \leq B \leq A,$ where $I$ is the identity operator.
It follows from formula (\ref{EqA.5}), that an explicit expression of $A$ is
given by
\begin{equation}
  (A \eta)(k)=E( U(\infty,\eta) u^{\infty}(k) | \mathcal{F}_{k})
 \label{EqA.6.1}
 \end{equation}
and that an explicit expression of $B$ is
given by
\begin{equation}
  (B \eta)(k)=E(( U(\infty,\eta) -E(U(\infty,\eta)))u^{\infty}(k) | \mathcal{F}_{k}),
 \label{EqA.6.2}
 \end{equation}
for $\eta \in \mathcal{H},$ where $0 \leq k \leq \bar{T}.$

Next result (Corollary 2.6 of \cite{T2}) solves the optimization  problem of this paragraph:
%%%%%%%%%%%%%%%%%%%%%%%%%%%%%%%%%%%%%%%
\begin{theorem}\label{ThA.4.1}
Let hypotheses ($H_{1}$), ($H_{2}$) and ($H_{3}$) be satisfied. If
$\mathcal{C}_{0}$ is non-empty, then optimization problem (\ref{EqA.0}) has
unique solution $\hat{\eta} \in \mathcal{C}_{0}.$
\end{theorem}
The solution $\hat{\eta}$ is given by a constructive approach in \cite{T2}.
In fact, in that reference a  Lagrangian formalism, an algorithm to invert the operators $A$ and $B$
and approximation methods for determining the multipliers are given.

\section{Proofs}
\label{Proof}
\label{Proof Main}
We first remind that if $(X_{n})_{n \geq 1}$ is a sequence of $\mathbb{R}^{k}$
valued random variables, uniformly bounded in $L^{2}$ (in its strong topology)
and d-convergent (i.e. in distributions)
to $X,$ then $(X_{n})_{n \geq 1}$ converges to $X$ in $L^{q},$ $1 \leq q <2.$
This follows using Skorohod's theorem (c.f. theorem 29.6 of \cite{Bill}),
uniform integrability and dominated convergence.
The following simple lemma will be used in the coming proofs. As earlier in this
paper let $\mathcal{H}$ be the Hilbert space of elements
$\eta \in \mathcal{E}(\mathbb{R}^{N}, \mathcal{A}),$ such that $\eta(t)=0$
for $t> \bar{T}.$
\begin{lemma}\label{LemM.1}
For $t \in \mathbb{N},$
we suppose that $(v(t,t'))_{t' \geq 0} \in \mathcal{E}(\mathbb{R}^{N}, \mathcal{A}),$
and that $\eta \in \mathcal{H}.$
Let $v(t,t')=0,$ for $t \geq t'$ and let
$\|E((v(k,t))^{2}| \mathcal{F}_{k}) \|_{L^{\infty}}$ $ < \infty,$ for $k < t.$
If $U'(t,\eta)=\sum_{0 \leq k \leq t} v(k,t) \cdot \eta(k),$ then 
\[E((U'(t,\eta))^{2}) \leq
(\sum_{0 \leq k \leq t}\|E((v(k,t))^{2}| \mathcal{F}_{k}) \|_{L^{\infty}})^{1/2}
 (\sum_{0 \leq k \leq t}E(|\eta|^{2})^{1/2},\]
\[|E(U'(t,\eta))| \leq
(\sum_{0 \leq k \leq t}\|E(v(k,t)| \mathcal{F}_{k}) \|^{2}_{L^{\infty}})^{1/2}
 (\sum_{0 \leq k \leq t}E(\eta)^{2})^{1/2}\]
and the map $\eta \mapsto U'(t,\eta)$ is d-continuous on bounded subsets of
$\mathcal{H}.$
\end{lemma}
Proof: The following calculus proves the first inequality:
$(E((U'(t,\eta))^{2}))^{1/2}$ $ \leq \sum_{0 \leq k \leq t}
    (E((v(k,t) \cdot \eta(k))^{2}))^{1/2}$
     $\leq \sum_{0 \leq k \leq t} (E(E(|v(k,t)|^{2}| \mathcal{F}_{k}) |\eta(k))|^{2}))^{1/2}$ \\
$\leq \sum_{0 \leq k \leq t} \|E(|v(k,t)|^{2}| \mathcal{F}_{k})\|_{L^{\infty}}^{1/2}
                   (E(|\eta(k))|^{2}))^{1/2}$ \\
$\leq (\sum_{0 \leq k \leq t} \|E(|v(k,t)|^{2}| \mathcal{F}_{k})\|_{L^{\infty}})^{1/2}
        (\sum_{0 \leq k \leq t} E(|\eta(k))|^{2}))^{1/2}.$ \\
Similarly the second inequality follows from
 $|E(U'(t,\eta))| \leq $ \\ $ \sum_{0 \leq k \leq t} |E(E(v(k,t)| \mathcal{F}_{k})
               \cdot \eta(k))|$
$\leq  \sum_{0 \leq k \leq t}  \|E(v(k,t)| \mathcal{F}_{k})\|_{L^{\infty}} |E(\eta(k))|$ \\
$\leq  (\sum_{0 \leq k \leq t}  \|E(v(k,t)| \mathcal{F}_{k})\|_{L^{\infty}}^{2})^{1/2}$
                          $(\sum_{0 \leq k \leq t}     |E(\eta(k))|^{2})^{1/2}.$

To prove the last statement, let $(\eta_{n})_{n \geq 1}$ be a sequence  in
$\mathcal{H},$ d-convergent to $\eta \in \mathcal{H}.$ It is enough to prove that
$f(U'(t,\eta_{n}))$ converges to $f(U'(t,\eta))$ for each continuous real function
$f$ with compact support (c.f. theorem 2.5.2 of \cite{Ito}). Let $f$ be a given
such function with compact support $K.$ For given $\epsilon > 0,$
there exists
a $C^{\infty}$ function $g$ be with compact support, being a subset of in $K,$
such that
$|f(x)-g(x)| \leq \epsilon /4,$ for $x \in \mathbb{R}.$ There exists $M \in \mathbb{R},$
such that $|g(x)-g(y)| \leq M |x-y|,$ for $x,y \in K.$ Let
$\Delta_{n}= |E(f(U'(t,\eta_{n}))-f(U'(t,\eta)))|.$ It follows that
$\Delta_{n} \leq |E(f(U'(t,\eta_{n}))-g(U'(t,\eta_{n})))|
    +|E(f(U'(t,\eta))-f(U'(t,\eta)))| +|E(g(U'(t,\eta_{n}))-g(U'(t,\eta)))|
      \leq \epsilon /2 +M E(|U'(t,\eta_{n})-U'(t,\eta)|).$
If $p^{-1} +q^{-1}=1,$ $1 \leq p < \infty,$ then
$ E(|U'(t,\eta)|)$ \\ $\leq
(\sum_{0 \leq k \leq t}  E(|v(k,t)|^{p})^{1/p}
(\sum_{0 \leq k \leq t} E|(\eta(k))|^{q})^{1/q},$
according to H\"older's inequality.
If  $Q_{q}=E(|v(k,t)|^{p})^{1/p},$ then
  $\Delta_{n} \leq \epsilon /2 +MQ_{q}E(|\eta_{n}-\eta(k)|^{q})^{1/q},$
$1 <q.$ Let $1 <q <2.$
Since $E(|\eta_{n}-\eta(k)|^{q})^{1/q}$ converges to zero for $q<2$
(which was pointed out in the beginning of this section), there exists $n_{0} \geq 1,$
such that $MQ_{q}E(|\eta_{n}-\eta(k)|^{q})^{1/q} \leq \epsilon /2,$ for
$n \geq n_{0}.$ Hence $\Delta_{n} \leq \epsilon,$ for $n \geq n_{0},$ which proves
the assertion. \\

\noindent
\textbf{Proof of Theorem \ref{Th1.2}}
Let $(\eta, \Vec{K}(0), \Vec{D}) \in \mathcal{C}_{c'}$  and $\theta=\xi+\eta.$
For $K(t)=\sum_{1 \leq j \leq \aleph}K^{(j)}(t),$ we have
$\Psi(t,K,\theta,0) \leq \sum_{0 \leq k \leq t} P(K(t) <0).$
Condition ($c_{4}'$) is satisfied. Since $E(K(t)) \geq 0,$ it follows that
\[P(K(t) <0) \leq P(|K(t)- E(K(t))| > E(K(t))).\]
Chebyshev's inequality gives
\[\Psi(t,K,\theta,0) \leq \sum_{0 \leq k \leq t} \frac{V(K(k))}{(E(K(k)))^{2}}.\]
It follows from ($c_{4}'$) that
$\Psi(t,K,\theta,0) \leq \sum_{0 \leq k \leq t} \epsilon' (k).$
This proves, according to the hypotheses of the theorem, that
$\Psi(t,K,\theta,0) \leq \epsilon (t),$ so constraint ($c_{4}$) is satisfied.

In the case of ($c_{7}$), since ($c_{7}'$) is satisfied, it follows similarly,
with $\theta^{(j)}=\xi^{(j)}+\eta^{(j)},$ that
\[\Psi(t,K^{(j)},\theta^{(j)},m^{(j)}) \leq \sum_{0 \leq k \leq t}
\frac{V(K^{(j)}(k)-m^{(j)}(k,\theta^{(j)}))}{(E(K^{(j)}(k)-m^{(j)}(k,\theta^{(j)})))^{2}}.\]
This inequality and ($c_{7}'$) give that
$\Psi(t,K^{(j)},\theta^{(j)},m^{(j)}) \leq \epsilon^{(j)} (t),$
which shows that also ($c_{7}$) is satisfied. Therefore ($c_{1}$)--($c_{8}$)
are satisfied, so
$(\eta, \Vec{K}(0), \Vec{D}) \in \mathcal{C}_{c}.$ This proves the theorem. \\

\noindent
\textbf{Proof of Theorem \ref{Th1.3}}
For $x \in \mathbb{R}^{N^{(j)}} \cup \{\tau_{f}\},$ let  $p_{j}(\tau_{f})=0$
and $\lambda_{j} ( \tau_{f})=1,$ and if $x \in \mathbb{R}^{N^{(j)}}$ let
$p_{j}(x)=x$  and $\lambda_{j} (x)=0.$
Let $p=(p_{1}, \ldots ,p_{\aleph})$ and $\lambda=(\lambda_{1}, \ldots ,\lambda_{\aleph}).$
Let $\mathcal{E}^{2}_{\bar{T}}(\mathbb{R}^{N},\mathcal{A})$ be the subset of
elements $\eta \in \mathcal{E}^{2}(\mathbb{R}^{N},\mathcal{A}),$ such that
$\eta(t) =0$ for $t>\bar{T}.$ By definition , if $\eta \in \tilde{\mathcal{P}}_{u,\bar{T}}$
then $p \circ \eta \in \mathcal{E}^{2}_{\bar{T}}(\mathbb{R}^{N},\mathcal{A}),$ where
$N=\sum_{1 \leq j \leq \aleph} N^{(j)}$ and
$\lambda \circ \eta \in \mathcal{E}^{2}_{\bar{T}}(\mathbb{R}^{\aleph},\mathcal{A}).$ We
introduce the Hilbert space
$H_{0}=\mathcal{E}^{2}_{\bar{T}}(\mathbb{R}^{N},\mathcal{A})
 \oplus \mathcal{E}_{\bar{T}}^{2}(\mathbb{R}^{\aleph},\mathcal{A})$ and the
Hilbert space $H_{1}$ of elements
$( \Vec{K}(0), \Vec{D}) \in \mathbb{R}^{\aleph}
            \oplus \mathcal{E}_{\bar{T}}^{2}(\mathbb{R}^{\aleph},\mathcal{A}),$
such that $\Vec{D}(0)=0.$ Let $H=H_{0} \oplus H_{1}.$

We shall first prove that the set $A$ of elements
$(p \circ \eta, \lambda \circ \eta, \Vec{K}(0), \Vec{D}) \in H$ satisfying
constraints ($c_{1}$), ($c_{2}$), ($c_{4}'$) and ($c_{7}'$) is d-compact.
Formula (\ref{Eq1.3}) in the case of $ U^{(j)},$  property ($d_{2}$) in
case of $D$ and the property of $m^{(j)}$ before formula  (\ref{Eq1.7}) show
that the constraints ($c_{1}$), ($c_{2}$), ($c_{4}'$) and ($c_{7}'$) are
satisfied for $\theta=\xi+\eta$ if and only if they are satisfied with $p \circ \theta$
instead of $\theta.$ Let $A'$ be the subset of $A$ such that $\lambda \circ \eta=0$
and let $\bar{p}(p \circ \eta, \lambda \circ \eta, \Vec{K}(0), \Vec{D})
           =(p \circ \eta, 0, \Vec{K}(0), \Vec{D}).$
Let $\theta$ be such that $p \circ \theta =\theta.$ We note that
$U^{(j)}(\infty,\theta^{(j)})=U^{(j)}(\bar{T}+T,\theta^{(j)}).$
According to the hypothesis of the theorem and according to formula (\ref{Eq1.10})
it follows that
\begin{equation}
\begin{split}
(V(K(\bar{T}+T,\theta)))^{1/2}&
     =(V(U(\infty,\theta)-\sum_{1 \leq k \leq \bar{T}+T}D(k,\theta)))^{1/2} \\
\geq (V(U(\infty,\theta))^{1/2}&-(V(\sum_{1 \leq k \leq \bar{T}+T}D(k,\theta))^{1/2} \geq
      (1-c)(V(U(\infty,\theta))^{1/2}.
\end{split} \notag
\end{equation}
Since $\theta=\xi+\eta,$ where $\xi \in \mathcal{P}_{r}$ and $\eta \in \mathcal{P}_{u, \bar{T}},$
it follows that
\begin{equation} \label{EqM.1}
(V(U(\infty,\eta))^{1/2} \leq (1-c)^{-1}(V(K(\bar{T}+T,\theta)))^{1/2}
             +(V(U(\infty,\xi))^{1/2}.
\end{equation}
The first term on the right hand side is uniformly bounded for $\eta \in A',$
according to constraint ($c_{4}'$). Since
$V(\xi^{(j)}_{i}(k)u_{i}^{(j)\infty}(k))=|(\xi^{(j)}_{i}(k)|^{2}V(u_{i}^{(j)\infty}(k)) < \infty,$
it follows that the second term also is uniformly bounded on $A'.$ This shows that
$V(U(\infty,\eta)$ is uniformly bounded for $\eta \in A'.$ By re-indexing the processes
$\eta^{(j)}_{i}$ and $u_{i}^{(j)\infty}$ and by using the hypotheses
($h_{1}$)-- ($h_{4}$), it then follows from Theorem \ref{ThA.4} that there exists
a closed ball $B'(0,R_{1}),$ with center $0$ and radius $R_{1} < \infty,$ in
$\mathcal{E}^{2}_{\bar{T}}(\mathbb{R}^{N},\mathcal{A}),$ such that
$\eta \in B'(0,R_{1})$ for all $\eta$ satisfying the constraint ($c_{4}'$).

We have, according to equation (\ref{Eq1.10})
\begin{equation} \notag
\begin{split}
(V(\sum_{1 \leq k \leq t}D^{(j)}(k)))^{1/2}& \leq
   (V(K^{(j)}(t) -m^{(j)}(t,\theta^{(j)})))^{1/2} \\
     &+(V(U^{(j)}(t,\theta^{(j)})))^{1/2} +(V(m^{(j)}(t,\theta^{(j)})))^{1/2}.
\end{split}
\end{equation}
The first term on the right hand side is uniformly bounded in $\eta^{(j)},$
for $\eta \in B'(0,R_{1}),$ according to constraint ($c_{7}'$), the second
according to Lemma \ref{LemM.1} and the third according to the hypotheses of the 
theorem. This shows that there is $R' < \infty$ such that
\begin{equation} \label{EqM.2}
V(D^{(j)}(t)) <R',
\end{equation}
for all $D^{(j)}(t)$ satisfying conditions ($c_{4}'$) and ($c_{7}'$).

According to constraint ($c_{7}'$), we have that
$E(K^{(j)}(t)) \geq \delta^{(j)}(t)K(0)+E(m^{(j)}(t,\xi^{(j)}+\eta^{(j)})).$
Since $m^{(j)}(t,\xi^{(j)}+\eta^{(j)})$ is uniformly in bounded (hypotheses
of the theorem)
in $\mathcal{E}^{2}(\mathbb{R},\mathcal{A}),$ for $\eta \in B'(0,R_{1}),$
it follows that $E(K^{(j)}(t))$ is bounded from below on $A'.$ In particular,
constraint ($c_{1}$) then gives that $|K^{(j)}(0)|$ is bounded on $A'.$
Formula (\ref{Eq1.10}) and the fact that $E(U^{(j)}(t,\xi^{(j)}+\eta^{(j)}))$
is uniformly bounded for $\eta \in B'(0,R_{1}),$ (c.f. Lemma \ref{LemM.1})
show that $E(D^{(j)}(t)),$ $1 \leq t$ are bounded from above. $|E(D(t,\xi+\eta))|$
is uniformly bounded for $\eta \in B'(0,R_{1}),$
since, by the hypotheses of the theorem, $D(t,\xi+\eta)$ is uniformly bounded
in $\mathcal{E}^{2}(\mathbb{R},\mathcal{A}),$ for $\eta \in B'(0,R_{1}).$ Since
$D=\sum_{1 \leq j \leq \aleph}D^{(j)},$ according to ($c_{2}$), it follows
that $E(D^{(j)}(t,\xi^{(j)}+\eta^{(j)}))$ is bounded on $A'.$ According to
inequality (\ref{EqM.2}) it follows then that there exists $R_{2}<\infty$ such that
\begin{equation} \label{EqM.3}
\|D^{(j)}(t)\| + |K^{(j)}(0)| <R_{2},
\end{equation}
for all $(\eta, 0, \Vec{K}(0), \Vec{D}) \in A'.$ This inequality and the fact
that $\eta \in B'(0,R_{1}),$  shows that $A'$ is a bounded subset of $H.$

Let $(x_{n})_{n \geq 1}$ be a sequence in $A',$ which converges in the d-topology
of $H$ to $x'=(\eta', 0, \Vec{K}'(0), \Vec{D}') \in H.$ The limit element then
satisfies the constraints ($c_{1}$) and ($c_{2}$). In fact, the case of ($c_{1}$)
is trivial. By the hypotheses of the theorem, it follows that
$x \mapsto \delta(x)=D(t,\xi+\eta) -\sum_{1 \leq j \leq \aleph}D^{(j)}(t),$
is d-continuous on $A'.$ We then have
$E(|\delta(x')|) \leq \liminf_{n \rightarrow \infty} E(|\delta(x_{n})|)=0$
(c.f. theorem 25.11 of \cite{Bill}), so $x'$ satisfies ($c_{2}$).

We note that the functions
\begin{equation} \label{EqM.4}
x \mapsto U^{(j)}(t,\xi^{(j)}+\eta^{(j)}), \quad x \mapsto K^{(j)}(t)
\end{equation}
and also the functions $x \mapsto E(U^{(j)}(t,\xi^{(j)}+\eta^{(j)})),$
$x \mapsto E(m^{(j)}(t,\xi^{(j)}+\eta^{(j)}))$ and
$x \mapsto E(K^{(j)}(t)-m^{(j)}(t,\xi^{(j)}+\eta^{(j)}))$
are d-continuous on $A'.$ In fact $A'$ is a bounded subset of $H,$ so according
to Lemma \ref{LemM.1}, the first function in (\ref{EqM.4}) is d-continuous with
$\mathcal{E}^{2}$ bounded image. According to formula (\ref{Eq1.10}), this is
then also the case for the second function in (\ref{EqM.4}). The boundedness
in $\mathcal{E}^{2}$ of the image of $A'$ under the function
$x \mapsto m^{(j)}(t,\xi^{(j)}+\eta^{(j)})$ (by hypotheses of the theorem) and
under the functions in (\ref{EqM.4}) and the d-continuity of these functions then
give the d-continuity of the above expected values (c.f. theorem 25.12 of \cite{Bill}).
The function $x \mapsto f(x)=K^{(j)}(t)-m^{(j)}(t,\xi^{(j)}+\eta^{(j)})$ is
d-continuous.
Since the sequence $(f(x_{n}))_{n \geq 1}$ converges in the d-topology, it follows
that $E(|f(x')|^{2}) \leq \liminf_{n \rightarrow \infty} E(|f(x_{n})|^{2})$
(c.f. theorem 25.11 of \cite{Bill}). Because of the already proved d-continuity
on $A'$ of $x \mapsto E(f(x)),$ it follows that
\[V(f(x'))=E((f(x'))^{2}) -(E(f(x')))^{2} \leq
             \liminf_{n \rightarrow \infty}V(f(x_{n})).\]
Constraint ($c_{7}'$) then gives that
$V(f(x')) \leq {\epsilon'}^{(j)}(t) (\delta^{j}(t)K(0))^{2},$
which proves that
\begin{equation} \label{EqM.5}
V(K^{'(j)}(t)-m^{(j)}(t,\xi^{(j)}+\eta^{'(j)})) 
              \leq {\epsilon'}^{(j)}(t) (\delta^{j}(t)K(0))^{2}.
\end{equation}
Once more, by the already proved proved d-continuity, it follows from ($c_{7}'$)
that
\begin{equation} \label{EqM.6}
E(K^{'(j)}(t)-m^{(j)}(t,\xi^{(j)}+\eta^{'(j)})) 
              \geq  \delta^{j}(t)K(0).
\end{equation}
Inequalities (\ref{EqM.5}) and (\ref{EqM.6}) show that the limit point $x'$
satisfies constraint ($c_{7}'$). Similarly it is proved that $x'$ also
satisfies constraint ($c_{4}'$). This proves that $A'$ is d-closed.

Let $(x_{n})_{n \geq 1}$ be a sequence in $A'.$
Then$(x_{n})_{n \geq 1}$ is tight, since $A'$ is a bounded subset of $H.$
In fact, if the bound is $R$ and if
$f(a)= \inf_{X \in A'} P(|X| \leq a),$ then
$f(a)=\inf_{X \in A'}(1-P(|X| > a)) \geq 1 -a^{-2}\sup_{X \in A'}E(|X|^{2})
            \geq 1-(R/a)^{2}.$
Hence $\lim_{a \rightarrow \infty}f(a)=1.$
This proves that $A'$ is sequentially d-compact, since $A'$ is closed.  It is
then d-compact (c.f. theorem 2.5.3 of \cite{Ito}).

Let $(y_{n})_{n \geq 1},$
$y_{n}= (p \circ \eta_{n},\lambda \circ \eta_{n}, \Vec{K}_{n}(0), \Vec{D}_{n}),$
be a sequence in $A$ and let $x_{n}=\bar{p}( y_{n}).$ Extracting a subsequence
we can suppose that  $(x_{n})_{n \geq 1}$ converges in the d-topology to $x \in A',$
since $A'$ is d-compact. Let $g_{n}=\lambda \circ \eta_{n}.$ Then
$g_{n} : \Omega \rightarrow \{0,1\}^{\aleph}$ defines a bounded sequence in
$\mathcal{E}^{2}(\mathbb{R}^{\aleph},\mathcal{A}).$ Extracting a subsequence
we can suppose that  $(g_{n})_{n \geq 1}$ converges in the d-topology to
$g \in \mathcal{E}^{2}(\mathbb{R}^{\aleph},\mathcal{A}).$ It follows that
$g : \Omega \rightarrow \{0,1\}^{\aleph}$ (a.e.) (since $g_{n}^{(j)}$ is a
d-convergent sequence of characteristic functions; $g_{n}^{(j)}=(g_{n}^{(j)})^{2}$).
This proves that a subsequence of $(y_{n})_{n \geq 1}$ converges to $y \in A,$
so $A$ is d-compact.

We shall next prove that the subset $B \subset A,$ of all elements in $A,$ which
satisfies conditions ($c_{3}$), ($c_{5}$), ($c_{6}$) and ($c_{8}$), is
d-compact. As earlier in this proof let
$x=(p \circ \eta,\lambda \circ \eta, \Vec{K}(0), \Vec{D}) \in A$ and let
$(y_{n})_{n \geq 1}$ be a d-convergent sequence in $A$ with limit
 $y=(p \circ \eta', \lambda \circ \eta', \Vec{K}'(0), \Vec{D}').$ $y \in A,$
since $A$ is d-compact. Suppose that $B$ is non-empty, otherwise there is
nothing to prove.

{\em Constraint} ($c_{3}$): Let $y_{n}$ satisfy ($c_{3}$) for $n \geq 1.$
Due to formula (\ref{Eq1.3}), it is enough to
consider the value of the constraint at $x_{n}=\bar{p}(y_{n}).$ It is proved,
as in the case of the d-continuity of functions in and following (\ref{EqM.4})
that both sides of the inequality in ($c_{3}$) are d-continuous functions on $A'.$
This proves that ($c_{3}$) is satisfied by $\bar{p}(y),$ since it is satisfied by
$\bar{p}(y_{n}),$ $n \geq 1.$ Therefore it is satisfied by $y.$

{\em Constraint} ($c_{5}$): Let $y_{n}$ satisfy ($c_{5}$) for $n \geq 1.$
According to the hypotheses of the theorem, the function
$x \mapsto f_{\alpha}(x)=F_{\alpha}(t,\eta, \Vec{K}(0), \Vec{D})$ is sequentially
lower semicontinuous in the d-topology. This gives that
$f_{\alpha}(y) \leq \liminf_{n \rightarrow \infty} f_{\alpha}(y_{n}).$ It
follows from ($c_{5}$) that
$f_{\alpha}(y) \leq C_{\alpha}(t,K(0),\xi).$ This proves that ($c_{5}$)
is satisfied by $y.$

{\em Constraint} ($c_{6}$): Let $y_{n}$ satisfy ($c_{6}$) for $n \geq 1.$
The functions, $\eta^{(j)} \mapsto p_{j} \circ \eta^{(j)}$ from
$\tilde{\mathcal{P}}_{u}^{(j)}$ to $\mathcal{P}_{u}^{(j)},$ and
$\eta^{(j)} \mapsto \lambda_{j} \circ \eta^{(j)}$ from
$\tilde{\mathcal{P}}_{u}^{(j)}$ to
$\mathcal{E}^{2}(\mathbb{R},\mathcal{A}),$ are d-continuous. According
to the hypotheses of the theorem,
$\eta \mapsto (\underline{c}^{(j)}_{i}((p \circ \eta))(t,\omega)$ is then lower
semi continuous (a.e.). This is then also the case for the function
$\eta \mapsto (b^{(j)}_{i}(\eta))(t,\omega)=(1-\lambda_{j} \circ \eta^{(j)})
             (\underline{c}^{(j)}_{i}((p \circ \eta))(t,\omega).$
The first condition of ($c_{6}$) is equivalent to
$(b^{(j)}_{i}(\eta))(t,\omega) \leq ((p_{j} \circ \eta^{(j)})(t))(\omega)
                 < \infty$ (a.e.).
Since the left hand side is lower semi d-continuous, it follows that
$(b^{(j)}_{i}(\eta'))(t,\omega) \leq \liminf_{n \rightarrow \infty}
         (b^{(j)}_{i}(\eta_{n}))(t,\omega).$
Since the member in the middle is d-continuous (a.e.), it then follows that
$(b^{(j)}_{i}(\eta'))(t,\omega) \leq \eta^{'(j)})(t))(\omega)$ (a.e.).
This inequality shows that ($c_{6}$) is satisfied by $y.$

{\em Constraint} ($c_{8}$): Let $y_{n}$ satisfy ($c_{8}$) for $n \geq 1.$
Let
\[s^{(j)}(t,x)=\min_{0 \leq k \leq t} (K^{(j)}(k)- m^{(j)}(k,\xi^{(j)}+\eta^{(j)})),\]
where $K^{(j)}(k)$ is evaluated at $x \in A.$ According to the hypotheses of the
theorem $x \mapsto s^{(j)}(\cdot,x)$ is d-continuous from $A$ to
$\mathcal{E}^{2}(\mathbb{R},\mathcal{A}).$ We shall prove by contradiction that
($c_{8}$) is satisfied by $y.$ We remind that, for $x \in A,$
$(\eta^{(j)}(t))(\omega) \in {\mathbb{R}}^{N^{(j)}}$ or $(\eta^{(j)}(t))(\omega)=\tau_{f}.$
Suppose that for some $t' \geq 0,$ $(s^{(j)}(t',y))(\omega) < 0$ for all
$\omega \in Q',$ where $Q'$ is a $\mathcal{F}_{t'}$-measurable set, with
$P(Q') >0.$ Suppose also that $(\eta^{'(j)}(t'+1))(\omega) \neq \tau_{f},$ for
all $\omega \in Q'',$ where $Q'' \subset Q'$ is a $\mathcal{F}_{t'}$-measurable
set and $P(Q'') >0.$ Let $h$ be the characteristic function of $Q''.$ Then
$(h,s^{(j)}(t',y))_{L^{2}} \leq -r <0,$ for some $r,$ and
$(h, (\lambda_{j} \circ \eta^{'(j)})(t'))_{L^{2}}=0.$

The sequences $(y_{n})_{n \geq 1}$ and $(s^{(j)}(t,y_{n}))_{n \geq 1}$ converge
in $L^{1}.$ In fact, $(y_{n})_{n \geq 1}$ is bounded in $H,$ so the $L^{1}$
convergence follows as noted in the beginning of this section.
It follows (as after formula (\ref{EqM.4})), that the image of $A$ under the
map $x \mapsto s^{(j)}(t,x)$ is bounded and that this map is d-continuous on
$A.$ Therefore also
$(s^{(j)}(t,y_{n}))_{n \geq 1}$ converges in $L^{1}.$ If $\epsilon >0,$ then
there exists $n_{0} \geq 1$ such that
\begin{equation} \label{EqM.7}
(h,s^{(j)}(t',y_{n}))_{L^{2}} \leq -r+\epsilon <0
\end{equation}
and $0 \leq (h, (\lambda_{j} \circ \eta^{(j)}_{n})(t'))_{L^{2}} \leq \epsilon,$
for $n \geq n_{0}$ (where we have used that $(y_{n})_{n \geq 1}$ satisfies
($c_{8}$)). Since $((\lambda_{j} \circ \eta^{(j)}_{n})(t))(\omega) \in \{0,1\},$
it follows that $(\eta^{(j)}_{n}(t'))(\omega)=\tau_{f}$ on a subset $\tilde{Q}_{n}$
of $Q'',$ with measure $P(\tilde{Q}_{n}) \leq \epsilon.$ Let $\tilde{h}_{n}$ be the
characteristic function of $\tilde{Q}_{n}.$ There exists $R,$ such that
$\| s^{(j)}(t,y_{n})\|_{L^{2}} \leq R,$ uniformly in $t$ and $n.$ We have
$|(\tilde{h}_{n},s^{(j)}(t',y_{n}))_{L^{2}}|
          \leq \|\tilde{h}_{n}\|_{L^{2}} \|s^{(j)}(t',y_{n})\|_{L^{2}}
         \leq (\epsilon)^{1/2} R.$
Since constraint ($c_{8}$) is satisfied by $y_{n},$ $n \geq 1,$ it follows that
$(\tilde{h}_{n},s^{(j)}(t',y_{n}))_{L^{2}} \leq 0$ and that
$(h- \tilde{h}_{n},s^{(j)}(t',y_{n}))_{L^{2}} \geq 0.$
We obtain that $(h,s^{(j)}(t',y_{n}))_{L^{2}} \geq (\epsilon)^{1/2} R,$
which is in contradiction with inequality (\ref{EqM.7}), if $\epsilon >0$ is
sufficiently small.

This shows that if $t \geq (t^{(j)}_{f})(\omega)+1,$ then
$(\eta^{'(j)}(t'))(\omega)= \tau_{f},$ for $\omega \in \Omega$ (a.e.).
This proves that $y$ satisfies the first part of the constraint ($c_{8}$).

That  $y$ satisfies the second part of the constraint ($c_{8}$)
(i.e. $(\eta^{'(j)}(t))(\omega) \in \mathbb{R}^{N^{(j)}},$ if 
$t \leq (t^{(j)}_{f})(\omega)$ and
$(K^{'(j)}(t-1)- m^{(j)}(t-1,\xi^{(j)}+\eta^{'(j)}))(\omega) \geq 0$
is proved similarly.
This ends the proof of constraint ($c_{8}$).

We have now proved that $y \in B,$ so $B$ is a sequentially d-closed subset of $A.$
Thus $B$ is a sequentially d-compact, since $B \subset A$ and $A$ is compact.
Using the boundedness of $B$ in $H,$ it then follows that $B$ is d-compact.

Similarly as after formulas (\ref{EqM.4}), it follows that
$\eta \mapsto E(U(t,\xi+\eta))$ is d-continuous on bounded subsets of
$\tilde{\mathcal{P}}_{u,\bar{T}}.$ The map $x \mapsto f(x)=E(U(t,\xi+\eta)),$
from$A,$ is then d-continuous. Since $B \subset A$ is d-compact, it follows
that the map $B \ni x \mapsto f(x)=E(U(t,\xi+\eta)),$ attains its supremum,
$f(\hat{x}),$ at a point $\hat{x} \in B.$ This proves the theorem. \\

\noindent
\textbf{Proof of Theorem \ref{Th1.4}}
This is a direct consequence of Theorem \ref{ThA.4.1}.

\bibliographystyle{amsplain}

\end{document}